\documentclass[a4paper,11pt]{article}
\usepackage{latexsym,amsfonts,amsmath}

\def\1{\mathbf{1}}
\def\inter{\mathop{\cap}}
\def\NN{\mathbb{N}}

\def\RR{\mathbb{R}}

\newcommand{\widebar}[1]{\overline{#1}}

\def\union{\mathop{\cup}}
\def\inter{\mathop{\cap}}

\def\ds{\displaystyle}

\def\nl{\mbox{} \newline }

\newcounter{hypot}

\newcounter{assump}

\begin{document}
\newtheorem{theorem}{Theorem}[section]
\newtheorem{proposition}[theorem]{Proposition}
\newtheorem{lemma}[theorem]{Lemma}
\newtheorem{corollary}[theorem]{Corollary}
\newtheorem{definition}[theorem]{Definition}
\newtheorem{remark}[theorem]{Remark}
\newtheorem{conjecture}[theorem]{Conjecture}
\newtheorem{assumption}[theorem]{Assumption}
\newtheorem{example}{Example}

\bibliographystyle{plain}

\title{A Linear Programming Formulation for Constrained Discounted Continuous Control for Piecewise Deterministic Markov Processes
\thanks{This work was partially supported by project USP/COFECUB and USP/MaCLinC.} }

\author{ \mbox{ }
O.L.V. Costa
\thanks{This author received financial support from CNPq (Brazilian National Research
Council), grant 301067/09-0.}
\thanks{Author to whom correspondence should be sent to.}\\
\\ \small Departamento de Engenharia de Telecomunica\c c\~oes e
Controle \\ \small Escola Polit\'ecnica da Universidade de S\~ao
Paulo \\ \small CEP: 05508 900-S\~ao Paulo, Brazil.
\\ \small phone: 55 11 30915771; fax: 55 11 30915718.
 \\ \small
e-mail: oswaldo@lac.usp.br
\\
\and
F. Dufour
\thanks{This author was supported by ARPEGE program of the French National Agency of Research (ANR),
project ``FAUTOCOES'', number ANR-09-SEGI-004.}
\\
\small Universit\'e Bordeaux \\
\small IMB, Institut Math\'ematiques de Bordeaux \\
\small INRIA Bordeaux Sud Ouest, Team: CQFD \\
\small \small 351 cours de la Liberation \\
\small 33405 Talence Cedex, France \\ \small e-mail :
dufour@math.u-bordeaux1.fr
}
\maketitle
\vspace*{-3mm}
\begin{abstract}
This papers deals with the constrained discounted control of piecewise
deterministic Markov process (PDMPs) in general Borel spaces.
The control variable acts on the jump rate and transition measure, and the goal is to minimize the total expected discounted cost,
composed of positive running and boundary costs, while satisfying some constraints also in this form.
The basic idea is, by using the special features of the PDMPs, to re-write the problem via an
embedded discrete-time Markov chain associated to the PDMP and
re-formulate the problem as an infinite dimensional linear programming (LP) problem, via the occupation measures associated to the
discrete-time process.
It is important to stress however that our new discrete-time problem is not in the same framework
of a general constrained discrete-time Markov Decision Process and, due to that, some conditions are required to get the equivalence between the continuous-time
problem and the LP formulation.
We provide in the sequel sufficient conditions for the solvability of the associated LP problem, based on a generalization of
Theorem 4.1 in \cite{dufour12}. In the Appendix we present the proof of this generalization which, we believe, is of interest on its own. The paper is concluded with some examples to illustrate the obtained results.
\end{abstract}
\begin{tabbing}
\small \hspace*{\parindent}  \= {\bf Keywords:}
\begin{minipage}[t]{13cm}
Markov Decision Processes, Continuous-Time, Discounted Cost, Constraints, General Borel Spaces,
Linear Programming
\end{minipage}
\\
\> {\bf AMS 2000 subject classification:} \= Primary 90C40
\\ \> \> Secondary 93E20
\end{tabbing}
\newpage
\section{Introduction}

Piecewise deterministic Markov processes (PDMPs) were introduced in \cite{davis84} and \cite{davis93} as a general family of continuous-time
non-diffusion stochastic models, suitable for formulating many optimization problems in queuing and inventory systems, maintenance-replacement models, and many other areas of engineering and operations research. PDMPs are determined by three local characteristics: the flow $\phi$, the jump rate $\lambda$, and the
transition measure $Q$. Starting from $x$, the motion of the process follows the flow $\phi(x,t)$ until the first jump time
$T_1$, which occurs either spontaneously in a Poisson-like fashion with rate $\lambda$ or when the flow $\phi(x,t)$ hits the boundary
of the state space. In either case the location of the process at the jump time $T_1$ is selected by the transition measure
$Q(\phi(x,T_1),.)$ and the motion restarts from this new point as before. As shown in \cite{davis93}, a suitable choice of the state space and the local
characteristics $\phi$, $\lambda$, and $Q$ provides stochastic models covering a great number of problems of engineering and operations research (see, for instance, \cite{davis93}, \cite{davis87}).

The objective of this work is to study the discounted continuous-time constrained optimal control problem of PDMPs by using the linear programming (LP) approach.
Roughly speaking the formulation of the control problem is as follows. After each jump time, a control is chosen from
a control set (which depends on the state variable) and will act on the jump rate $\lambda$ and transition
measure $Q$ until the next jump time. The control variable will have two components, one that will parametrize a function that will
regulate the jump rate and transition measure before the flow hits the boundary, and the other component that will act on the
transition measure at the boundary (see Remark \ref{remarkA} below). The goal is to minimize the total expected discounted cost, which is composed of a
running cost and a boundary cost (added to the total cost each time the PDMP touches the boundary). Both costs are assumed to be positive but
not necessarily bounded. The constraints are also in the form of total expected discounted cost, again composed of
positive running and boundary costs, not necessarily bounded. The state and control spaces are assumed to be
general Borel spaces.

The linear programming technique  has proved to be a very efficient method for solving continuous-time Markov Decision Processes (MDPs) with constraints.
We do not attempt to present an exhaustive panorama on this topic, but instead we refer the interested reader to
\cite{Altamn,Borkar-02,GUO-2011,GUO-2003,Hernandez-Lerma-2000,Piunovskiy-97} and the
references therein for detailed expositions on this technique in the context of continuous-time controlled Markov processes.

Contrary to continuous-time constrained MDPs, it is important to emphasize that constrained optimal control problems of PDMPs
have received less attention.
An attempt in this direction is presented in \cite{goreac12} where the authors study a control problem for a special class of PDMPs (with no boudary) by using an LP technique.


In this paper we aim at tracing a parallel with the theory developed for discounted discrete-time constrained MDP in general Borel spaces in order to get our results.  We follow a similar approach as the one used for studying continuous-time MDPs without constraints, which consists of reducing
the original continuous-time control problem into a semi-Markov or discrete-time MDP.
For more details about such equivalence results, the reader may consult the recent survey \cite{feinberg12} and the manuscript \cite{feinberg04}, and references therein.
It is important to underline that such results cannot be directly applied in our case. Indeed, PDMPs are not piecewise constant processes and moreover, and more importantly, PDMPs have \textit{deterministic jumps} in the sense that, roughly speaking, the process necessarily jumps when it hits the boundary. As a consequence, the inter-arrival jumping times are not exponentially distributed and the compensator of the process is not absolutely continuous as for continuous-time MDPs.
Regarding PDMPs, the idea of reducing the original continuous-time control problem into a semi-Markov or discrete-time MDP
was developed in \cite{davis86} by reformulating the optimal control problem of a PDMP for a discounted cost as an equivalent discrete-time MDP in which the stages are the jump times $T_{n}$ of the PDMP.
In this paper we follow similar steps by re-writing the discounted
continuous-time control constrained PDMP as a constrained discrete-time problem, in function of the post-jump location and control action,
and with the stages being the jump times $T_{n}$ of the PDMP. By doing this we can trace a parallel with the general theory for discrete-time MDPs in Borel spaces. 
As usual in these problems, we define the space of occupation measures associated to the discrete-time problem and, from this, we re-formulate the original problem as an infinite dimensional linear programming  problem.
It is important to stress however that our new discrete-time problem is given in terms of an MDP with an expected total cost criterion, in particular there is no discount factor multiplying the stochastic kernel.
Due to that an expected growth condition is introduced to derive the equivalence between the continuous-time
problem and the LP formulation. Conditions for the solvability of the LP are also provided, based on Theorem 4.1 in \cite{dufour12}.
Note however that in Theorem 4.1 in \cite{dufour12} it is supposed that the control space $A$ is fixed and doesn't depend on the state variable,
while in the present paper the control space is dependent on the state variable. Consequently, Theorem 4.1 in \cite{dufour12}
is generalized in the appendix to an expected total cost discrete-time constrained MDP with the control space in the form $A(x)$, that is, depending on the state variable $x$.
We believe that this result, presented in the appendix in the general context of discrete-time constrained MDPs, is interesting on its own and could also be applied in other discrete-time constrained MDPs problems.

The paper is organized as follows. In section \ref{prob-form} we introduce the notation, some definitions and the problem formulation,
while in section \ref{MainOp} we present some important operators that will be needed to get our results,
some of the main assumptions and some auxiliary results. In section \ref{AuxRes} we present some auxiliary results establishing a connection between the
continuous-time problem and the discrete-time formulation, as well as the occupation measures that will be considered for the LP formulation.
In section \ref{sec6} we introduce the LP formulation and the first main result of the paper, Theorem \ref{consistency},  which establishes the equivalence between the original problem and the LP problem. In section \ref{Solvability} we present the second main result of the paper,
Theorem \ref{theoSol}, which provides conditions for the solvability of the LP problem. The paper is concluded with section \ref{Exam}, where we first
present the finite dimensional LP formulation of the problem for the case in which the control space and post-jump location set are finite, and in the sequel
we illustrate the application of our formulation and conditions for a capacity expansion problem. In the Appendix it is presented a generalization
of Theorem 4.1 in \cite{dufour12}, as explained before.

\section{Notation and Problem Formulation}
\label{prob-form}
\subsection{Notation}
In this section we present some standard notation, basic definitions and the main assumptions related to the motion of a PDMP $\{X(t)\}$,
and the control problem we will consider throughout the paper. For further details and properties the reader is referred to \cite{davis93}. The following notation will be used in this paper:
$\NN$ is the set of natural numbers, $\NN^{*}=\NN-\{0\}$, $\NN_p=\{0,\ldots,p\}$ and $\NN_p^{*}=\NN_p-\{0\}$. $\RR$ denotes the set of
real numbers, $\RR_+$ the set of positive real numbers, $\widebar{\RR}_{+}=\RR_{+}\union \{+\infty\}$ and
$\RR^d$ the $d$-dimensional euclidian space. 
For $X$ a metric space, we denote by $\mathcal{B}(X)$ the $\sigma$-algebra generated by the open sets of $X$.
$\mathcal{M}(X)$ (respectively, $\mathcal{M}(X)_{+}$, $\mathcal{P}(X)$)
denotes the set of all finite signed (respectively, positive, probability) measures on
$(X,\mathcal{B}(X))$. Let $X$ and $Y$ be metric spaces. The set of all Borel measurable (respectively, bounded) functions from $X$ into $Y$ is denoted by $\mathbb{M}(X;Y)$ (respectively, $\mathbb{B}(X;Y)$).
Moreover, for notational simplicity $\mathbb{M}(X)$ (respectively, $\mathbb{B}(X)$, $\mathbb{M}(X)_{+}$, $\mathbb{B}(X)_{+}$)
denotes $\mathbb{M}(X;\RR)$ (respectively, $\mathbb{B}(X;\RR)$, $\mathbb{M}(X;\RR_{+})$, $\mathbb{B}(X;\RR_{+})$).
$\mathbb{C}_b(X)$ denotes the set of continuous bounded functions from $X$ to $\RR$.
For $g\in \mathbb{M}(X)$ with $g(x)>0$ for all $x\in X$, $\mathbb{B}_{g}(X)$ is the set of functions  $v\in \mathbb{M}(X)$ such that $\ds ||v||_{g}=\sup_{x\in X} \frac{|v(x)|}{g(x)}< +\infty$, and $\mathcal{M}_{g}(X)$ the set of finite signed measures $\mu \in  \mathcal{M}(X)$ such that $\|\mu\|_g = \int g d|\mu| < \infty$. As before we denote by $\mathbb{B}_{g}(X)_{+} = \{v\in \mathbb{B}_{g}(X): v \geq 0 \}$ and $\mathcal{M}_{g}(X)_{+}=\{\mu \in \mathcal{M}_{g}(X): \mu \geq 0\}$, and $\mathcal{P}_{g}(X) = \mathcal{P}(X)\cap \mathcal{M}_{g}(X)_{+}$.
The set of positive measures (not necessarily finite) on the metric space $X$ is denoted by $\mathfrak{M}(X)_{+}$.

If $X$ and $Y$ are Borel spaces, the stochastic kernel on $X$ given $Y$ is a function $F(.;.)$ such that $F(y;.)$ is a probability measure on $X$ for each fixed $y\in Y$, and $F(.;B)$ is a measurable function on $Y$ for each fixed $B\in \mathcal{B}(X)$. The family of all stochastic kernels on $X$ given $Y$ is denoted by $\mathcal{P}(X|Y)$.

Finally we define $I_{A}$ as the indicator function of the set $A$, that is, $I_{A}(x)=1$ if $x\in A$ and $I_{A}(x)=0$, otherwise, and
$ \delta_{x}$ as the Dirac measure centered on a fixed point $x\in E$.

\subsection{Problem Formulation}
\label{prob-form-PDMP}
For the definition of the state space of the PDMP we will consider for notational simplicity that $E$ is an open subset of $\RR^n$ with boundary $\partial E$ and closure $\widebar{E}$.
This definition could be easily
generalized to include some boundary points and countable union of sets as in \cite[section 24]{davis93}.
In what follows the set $\mathbb{U}$ of control actions is a Borel space. For each $x\in \widebar{E}$, we define the subsets $\mathbb{U}(x)$
of $\mathbb{U}$ as the set of feasible control actions that can be taken when the state process is in $x\in \widebar{E}$.

\bigskip
Let $\widetilde{\mathbb{U}}$ be another Borel space and consider a function $\ell \in  \mathbb{M}(E\times \mathbb{U}\times \mathbb{R}_{+},\widetilde{\mathbb{U}})$.
We introduce next some data that will be used to define the controlled PDMP.
\begin{itemize}
\item The flow $\phi(x,t)$ is a function $\phi:  \mathbb{R}^{n}\times \RR_{+} \longrightarrow \mathbb{R}^{n}$ continuous in $(x,t)$ and such that
$\phi(x,t+s) = \phi(\phi(x,t),s)$.
\item For each $x\in E$, the time the flow takes to reach the boundary starting from $x$ is defined as
$$t_{*}(x)\doteq \inf \{t>0:\phi(x,t)\in \partial E \}.$$
For $x\in E$ such that
$t_{*}(x)=\infty$ (that is, the flow starting from $x$ never touches the boundary), we set $\phi(x,t_{*}(x))=\Delta$, where $\Delta$ is a fixed point in $\partial E$.
\item The jump rate $\lambda \in \mathbb{M}(\widebar{E}\times\widetilde{\mathbb{U}})_{+}$.
\item The transition measure $Q$ which is a stochastic kernel $Q \in \mathcal{P}\big(E|(E\times \widetilde{\mathbb{U}})\union (\partial E\times \mathbb{U}) \big)$.
To avoid jumps to the same point, we assume that $Q(x,\widetilde{a};\{x\})=0$ for any $x\in E$, $\widetilde{a}\in \widetilde{\mathbb{U}}$.
\end{itemize}

The following assumption, based on the standard theory of MDPs (see \cite{hernandez96}),
will be made throughout the paper.

\begin{assumption}\label{Hyp2a}
The set
\begin{equation*}
\mathcal{K}=\bigl\{ (x,(a,a_\partial)): x\in E, a \in \mathbb{U}(x),\,\,a_\partial \in \mathbb{U}(\phi(x,t_{*}(x))) \bigr\} \subset E \times \mathbb{U} \times \mathbb{U}
\end{equation*}
is a Borel subset of  $E \times \mathbb{U} \times \mathbb{U}$.
\end{assumption}

\begin{remark}\label{remarkA} We will consider in this paper a simplified version of the controlled PDMP with respect to the one adopted in \cite{book-Springer}.
The idea is that after a jump from a point $x\in E$ two actions will be chosen, one from  $a\in \mathbb{U}(x)$ and another one from $a_\partial \in \mathbb{U}(\phi(x,t_*(x)))$.
Action $a_\partial$ will regulate the transition measure at the frontier point $\phi(x,t_*(x))$, while action $a$ will parametrize
a function $\ell(x,a,t)$, with $0 \leq t < t_*(x)$,  which will regulate the jump rate and transition measure of the PDMP until the next jump time
In \cite{book-Springer} it was supposed that the controller could freely choose a function to regulate the jump rate and transition measure of the PDMP between jump times,
instead of being restricted to the function $\ell(x,a,t)$. Therefore in the model considered in this paper
the decision is taken only after a jump time, and the control that will be applied to $\lambda$ and $Q$ will be restricted to a pre-defined function $\ell(x,a,t)$
for $0\leq t < t_*(x)$.
\end{remark}

We make the following definitions. For $x\in E$, define
\begin{align*}
\mathbb{S}(x)& = \mathbb{U}(x)\times \mathbb{U}(\phi(x,t_{*}(x))),\,\,\,\mathbb{S}^{r}(x) = \mathcal{P}\bigl(\mathbb{S}(x)\bigr),\\
\mathbb{S} & = \mathbb{U} \times \mathbb{U},\,\,\,\mathbb{S}^{r}  =  \mathcal{P}\bigl(\mathbb{S}\bigr).
\end{align*}

\begin{definition}
\label{Def-Mes-Selec}
\index{Selector}
The following set of measurable selectors will be considered:
\begin{eqnarray*}
\mathcal{S}_{\mathbb{S}^{r}} &=& \bigl\{ \mu\in \mathbb{M}(E;\mathbb{S}^{r}); \text{ for all }x\in E,\,\mu(x) \in \mathbb{S}^{r}(x)\bigr\}.
\end{eqnarray*}
\end{definition}

\begin{definition}
\label{DefAdcontStra-U}
An admissible randomized control strategy $U$ is defined as $U=\{\eta_k;k\in \NN\}$ such that for each $k\in \NN$, $\eta_k \in\mathcal{S}_{\mathbb{S}^{r}}$.
The class of admissible randomized control strategies will be denoted by $\mathcal{U}$.
\end{definition}

\begin{definition}
We say that $U$ is an admissible randomized stationary control strategy if $U=\{\eta_k;k\in \NN\}$ is such that for all $k\in\NN$, $\eta_k =\varphi$ for some
$\varphi \in\mathcal{S}_{\mathbb{S}^{r}}$. We denote the set of admissible stationary control strategies by $\mathcal{U}_s$, and write in this case $U=U^\varphi$.
\end{definition}

Given an admissible randomized control strategy $U=\{\eta_k\}$ we present next the definition of the controlled piecewise deterministic Markov process.
Consider the state space $\widehat{E}=E\times \mathcal{K} \times \RR_{+} \times \NN $.
Let us introduce the following parameters for $\hat{x}=(x,z,a,a_\partial,s,n)\in \widehat{E}$, where $(z,a,a_\partial)\in \mathcal{K}$:
\begin{itemize}
\label{flowflow}
\item the flow $\widehat{\phi}(\hat{x},t) = (\phi(x,t),z,a,a_\partial,s+t,n)$,
\item the jump rate $\widehat{\lambda}^{U}(\hat{x})=\lambda(x,\ell(z,a,s))$,
\item the transition measure
\[
\widehat{Q}^{U}(\hat{x};A\times B \times C \times \{0\}\times \{n+1\}) =
\begin{cases}
\int_{A\inter B} \eta_{n+1}(y;C) Q(x,\ell(z,a,s); dy ) & \text{if } x\in E,\\
\int_{A\inter B} \eta_{n+1}(y;C) Q(x,a_\partial; dy ) & \text{if } x\in \partial E,
\end{cases}
\]
for $A,B \in \mathcal{B}(E)$, and $C \in \mathcal{B}(\mathbb{S})$.
\end{itemize}
From \cite[section 25]{davis93}, it can be shown that for any admissible control strategy $U=\{\eta_k;k\in \NN\}$
there exists a filtered probability space $(\Omega,\mathcal{F},\{ \mathcal{F}_{t} \}, \{ P^{U}_{\hat{x}} \}_{\hat{x}\in \widehat{E}})$
such that the PDMP $\{\widehat{X}^{U}(t)\}$ with local characteristics $(\widehat{\phi},\widehat{\lambda}^{U},\widehat{Q}^{U})$ may be constructed as follows.
For notational simplicity the probability $P^{U}_{\hat{x}_{0}}$ will be denoted by $P^{U}_{(x,\widehat{a},k)}$ for $\hat{x}_{0}=(x,x,\widehat{a},0,k)\in \widehat{E}$, $\widehat{a}=(a,a_\partial)$.
Moreover, $E^{U}_{\hat{x}_{0}}$ denotes the expectation under the probability $P^{U}_{\hat{x}_{0}}$ and
$E^{U}_{\hat{x}_{0}}$ will be denoted by $E^{U}_{(x,\widehat{a},k)}$ for $\hat{x}_{0}=(x,x,\widehat{a},0,k)\in \widehat{E}$. Take a random variable $T_1$ such that
\begin{equation}\label{distri-T}
P^{U}_{(x,\widehat{a},k)}(T_1>t) \doteq
\begin{cases}
e^{-\Lambda^{U}(x,a,k,t)} & \text{for } t<t_{*}(x), \\
0 & \text{for } t\geq t_{*}(x),
\end{cases}
\end{equation}
where for $x\in E$ and $t\in [0,t_*(x) \mathclose[$,
$\Lambda^{U}(x,a,k,t) \doteq \int_0^t\lambda(\phi(x,s),\ell(x,a,s))ds$.
If $T_1$ is equal to infinity, then for $t\in \RR_+$, $\widehat{X}^{U}(t)= \bigl(\phi(x,t),x,\widehat{a},t,k\bigr)$. Otherwise
select independently an $\widehat{E}$-valued random variable
(labeled $\widehat{X}^{U}_{1}$) having
distribution\[
P^{U}_{(x,\widehat{a},k)}(\widehat{X}^{U}_{1} \in A \times B \times C \times \{0\} \times \{k+1\} |\sigma\{T_1\})
\]
(where $\sigma\{T_1\}$ is the $\sigma$-field generated by $T_1$) defined
by
\begin{equation}\label{distri-Q}
\begin{cases}
\int_{A\inter B}\eta_{k+1}(y;C) Q(\phi(x,T_{1}),\ell(x,a,T_{1}); dy ) & \text{if } \phi(x,T_{1})\in E, \\
\int_{A\inter B}\eta_{k+1}(y;C) Q(\phi(x,T_{1}),a_{\partial}; dy) & \text{if } \phi(x,T_{1}) \in \partial E.
\end{cases}
\end{equation}
The trajectory of $\{\widehat{X}^{U}(t)\}$ starting from $(x,x,\widehat{a},0,k)$, for $t\leq T_1$, is given
by\[
\widehat{X}^{U}(t) \doteq
\begin{cases}
\bigl(\phi(x,t),x,\widehat{a},t,k\bigr) &\text{for } t<T_1, \\
\widehat{X}^{U}_1 &\text{for } t=T_1.
\end{cases}
\]
Starting from $\widehat{X}^{U}(T_1)=\widehat{X}^{U}_1$, we now
select the next interjump time $T_2-T_1$ and post-jump location
$\widehat{X}^{U}(T_2)=\widehat{X}^{U}_2$ in a similar way.
The sequence of jump times of the PDMP is denoted by $(T_{n})_{n\in \NN}$.
Let us define the components of the PDMP $\{\widehat{X}^{U}(t)\}$ by
\begin{eqnarray}
\widehat{X}^{U}(t)=\bigl(X(t),Z(t),A(t),\tau(t),N(t)\bigr).
\label{defXU}
\end{eqnarray}
From the previous construction of the PDMP $\{\widehat{X}^{U}(t)\}$, it is easy to see that $X(t)$ corresponds to the trajectory of the
system, $Z(t)$ is the value of $X(t)$ at the last jump time before $t$, $A(t)$ the control action that will be applied from
the last jump time before $t$ until the next jump time and at the boundary, $\tau(t)$ is time elapsed between the last jump and time $t$,
and $N(t)$ is the number of jumps of the process  $\{X(t)\}$ at time $t$. It will be convenient to define $T_0=0$, $Z_0=X(0)$ and
\begin{enumerate}
\item [a)] $Z_i = X(T_i) \in E \rightarrow$ the post jump location after the $i^{th}$ jump,
\item [b)] $\Theta_i = (\theta_i,\theta_{i,\partial}) = A(T_i) \in \mathbb{S}(Z_i)\rightarrow$ the control action that will be applied from the $i^{th}$ jump until the $({i+1})^{th}$ jump.
\end{enumerate}

\begin{definition}
For $\nu_0\in \mathcal{P}(E)$, we define $P^{U}_{\nu_0}(D)$ for any $D\in \mathcal{B}(\widehat{E})$ as
\begin{align}\label{expect-nu0}
P^{U}_{\nu_0}(D) &= \int_{\mathcal{K}} P^{U}_{(x,\widehat{a},0)}(D) \eta_0(x;d\widehat{a}) \nu_0(dx) = \int_{E} \int_{\mathbb{S}(x)} P^{U}_{(x,\widehat{a},0)}(D) \eta_0(x;d\widehat{a}) \nu_0(dx).
\end{align}
We denote by $E^{U}_{\nu_0}(.)$ the expectation under the probability $P^{U}_{\nu_0}(.)$. We just write $E^{U}_{x_0}(.)$ and $P^{U}_{x_0}(.)$
for the case in which $\nu_0$ is the dirac measure over $x_0\in E$.
\end{definition}

The cost and restrictions of our control problem will contain two terms, a running cost $f_i$ and a boundary cost $r_i$, $i=0,1,\ldots,n$, satisfying the following properties.

\begin{assumption}\label{Hyp6a}
$f_i\in \mathbb{M}(E\times\mathbb{U})_{+}$ and $r_i\in \mathbb{M}(\partial E\times\mathbb{U})_{+}$, $i=0,1,\ldots,n$.
\end{assumption}

Define for $\alpha > 0$, $t\in \mathbb{R}_+$, and $U\in \mathcal{U}$, the finite-horizon $\alpha$-discounted cost functions:
\begin{align*}
\mathbf{J}_i^\alpha(U,t) = & {\int_{0}^{t} e^{-\alpha s}f_i\bigl(X(s), \theta_{N(s)}\bigr) ds} + {\int_{]0,t]}e^{-\alpha s} r_i\bigl(X(s-),\theta_{N(s-),\partial} \bigr) dp^{*}(s),}
\end{align*}
where $p^{*}(t) = \sum_{i=1}^{\infty} I_{\{T_{i}\leq t\}} I_{\{ X(T_{i}-) \in \partial E\}}$ counts the number of times the
process hits the boundary up to time $t$.
The sequence $\{T_{k}\}_{k\in \NN}$ is increasing and let us denote by $T_{\infty}=\lim_{k\rightarrow\infty} T_{k}\in \widebar{\RR}_{+}$.

For any $\nu_0\in \mathcal{P}(E)$ the infinite-horizon expected $\alpha$-discounted costs are given by
\begin{equation}
\mathcal{D}_i^{\alpha}(U,\nu_0) = E^{U}_{\nu_0}[\mathbf{J}_i^\alpha(U,T_{\infty})], \quad i\in \NN_{n}
\label{al-disc0}
\end{equation}
and the constrained $\alpha$-discount value function is defined by
\begin{equation}
\mathcal{J}_{\mathcal{D}}^{\alpha}(\nu_0) =
\inf\{\mathcal{D}_0^{\alpha}(U,\nu_0); \mathcal{D}_i^{\alpha}(U,\nu_0)\leq d_i, \quad i\in \NN_{n}^{*}, U\in\mathcal{U}\},
\label{al-disc1}
\end{equation}
where $(d_{i})_{i\in \NN_{n}^{*}}$ are the constraint limits. We need the following assumption to avoid infinite costs for the discounted case.

\begin{assumption}\label{Hypdis}
$\mathcal{J}_{\mathcal{D}}^{\alpha}(\nu_0)<\infty$.
\end{assumption}

We conclude this section with the following definition and result that will be used in section \ref{sec6}:
\begin{definition}
We define $\mathbb{M}^{ac}(E)$ as the space of functions absolutely continuous along the flow with limit towards the
boundary:\begin{align*}
\mathbb{M}^{ac}&(E)  =  \Bigl\{  g\in {\mathbb{M}(E);}\text{ }g(\phi(x,t)): [0,t_{*}(x)) \mapsto \RR \text{ is absolutely continuous for} \\
&\text{each }  x\in E \text{ and when } t_{*}(x)< \infty \text{ the limit } {\lim_{t\rightarrow t_{*}(x)} g(\phi(x,t))} \text{ exists in } \RR \Bigr\}.
\end{align*}
For $g\in \mathbb{M}^{ac}(E)$ and $z\in \partial E$ for which there exists $x\in E$ such that $z=\phi(x,t_{*}(x))$, where $t_{*}(x)< \infty$,
we define $g(z) = \lim_{t\rightarrow t_{*}(x)} g(\phi(x,t))$ (note that the limit exists by assumption).
\end{definition}
For any function $g\in \mathbb{M}^{ac}(E)$ we introduce in the next lemma the function $\mathcal{X}g$.
$\mathcal{X}$ can be seen as the vector field associated to the flow $\phi$.
The proof of this lemma can be found in \cite{book-Springer}.
\begin{lemma}
\label{Iexist} Assume that $w\in \mathbb{M}^{ac}(E)$. Then there exists a function $\mathcal{X}w$ in $\mathbb{M}(E)$ such that
for all $x\in E$, and $t\in [0,t_{*}(x))$
\begin{eqnarray}
w(\phi(x,t))-w(x) & = & \int_{0}^{t} \mathcal{X}w(\phi(x,s)) ds.
\label{derixt}
\end{eqnarray}
\end{lemma}

\section{Main Operators}\label{MainOp}
In this section we present some important operators associated to the constrained optimal control problem posed in (\ref{al-disc0}).
We will need the following assumption.

\begin{assumption}\label{Hyp4a}
There exists  $\widebar{\lambda}\in \mathbb{M}(\widebar{E})$, $\underline{\lambda} \in\mathbb{M}(\widebar{E})_{+}$
and $K_{\lambda} \in \RR_{+}$ such that, for any $(x,\widetilde{a})\in E\times \widetilde{\mathbb{U}}$,
\begin{enumerate}
\item[{\rm(a)}]
$\lambda(x,\widetilde{a})\leq \widebar{\lambda}(x)$, and
for $t\in [0,t_{*}(x))$, $ \int_{0}^{t} \widebar{\lambda}(\phi(x,s))  ds < \infty$, and if $t_{*}(x)< \infty$,
then $$\int_{0}^{t_{*}(x)} \widebar{\lambda}(\phi(x,s))  ds < \infty.$$
\item[{\rm(b)}] $\lambda(x,\widetilde{a}) \geq \underline{\lambda}(x)>0$ and
$\int_0^{t_{*}(x)} e^{-\int_0^t \underline{\lambda}(\phi(x,s))ds}dt \leq K_{\lambda}$.
\end{enumerate}
\end{assumption}

For any $h\in \mathbb{M}(\mathcal{K})_{+}$ and $\eta\in \mathcal{S}_{\mathbb{S}^{r}}$, we introduce the following notation:
\begin{align}
\label{notation-hvarphi}
h(x,\eta) & \doteq  \int_{\mathbb{S}(x)}h(x,\widehat{a}) \eta(x,d\widehat{a}).
\end{align}

We define for $x\in E$, $0\leq t < t_{*}(x)$, $\widehat{a}=(a,a_\partial) \in \mathbb{S}(x)$ and $A\in \mathcal{B}(E)$:
\begin{align}
\Lambda^{a}(x,t)   \doteq& {\int_{0}^{t} \lambda(\phi(x,s),\ell(x,a,s)) ds}, \label{Lambda-mu}\\
\lambda QI_{A}(\phi(x,t),\ell(x,a,t)) \doteq {} & \nonumber\lambda(\phi(x,t),\ell(x,a,t)) QI_{A}(\phi(x,t),\ell(x,a,t))
\end{align}
Let us introduce the kernel $G$ on $E$ given $\mathcal{K}$ as follows:
\begin{align}
\label{DefGr}
G(x,\widehat{a};A)   \doteq  &   {\int_0^{t_{*}(x)}e^{-\alpha s - \Lambda^{a}(x,s)}\lambda QI_{A}(\phi(x,s),\ell(x,a,s)) ds} \nonumber \\
&  {+}\> e^{-\alpha t_{*}(x) -\Lambda^{a}(x,t_{*}(x))} Q(\phi(x,t_{*}(x)),a_{\partial};A),
\end{align}
Clearly we have that $G(x,\widehat{a};A)\leq 1$ for any $(x,\widehat{a})\in \mathcal{K}$ and $A\in \mathcal{B}(E)$. Now introduce the operator
$L$ (respectively, $H$) defined on $\mathbb{M}(E\times \mathbb{U})_{+}$ (respectively, $\mathbb{M}(\partial E\times \mathbb{U})_{+}$) with values in
$\mathbb{M}(E\times \mathcal{K};\widebar{\RR}_{+})$ (respectively, $\mathbb{M}(E\times \mathcal{K};\RR_{+})$) as follows:
\begin{align}
\label{DefLr}
Lv(x,\widehat{a}) & \doteq   \int_0^{t_{*}(x)}e^{-\alpha s-\Lambda^{a}(x,s)} v(\phi(x,s),a) ds, \\
\label{DefHr} Hw(x,\widehat{a}) & \doteq  e^{-\alpha t_{*}(x)-\Lambda^{a}(x,t_{*}(x))} w(\phi(x,t_{*}(x)),a_{\partial}),
\end{align}
for $v\in \mathbb{M}(E\times \mathbb{U})_{+}$, $w\in \mathbb{M}(\partial E\times \mathbb{U})_{+}$. For $h\in \mathbb{M}(E)$ (respectively, $v\in \mathbb{M}(E\times \mathbb{U})$),  $Gh(x,\widehat{a})=Gh^{+}(x,\widehat{a})-Gh^{-}(x,\widehat{a})$
(respectively, $Lv(x,\widehat{a})=Lv^{+}(x,\widehat{a})-Lv^{-}(x,\widehat{a})$) provided the difference has a meaning.
By a slight abuse of notation, we also write for $h\in \mathbb{M}(E)_{+}$
$$Lh(x,\widehat{a})\doteq\int_0^{t_{*}(x)}e^{-\alpha s - \Lambda^{a}(x,s)} h(\phi(x,s))ds$$
and for $g\in \mathbb{M}(E)$, $Lg(x,\widehat{a})=Lg^{+}(x,\widehat{a})-Lg^{-}(x,\widehat{a})$ provided the difference has a meaning.

\begin{remark}\label{vac}\rm
From a) of Assumption \ref{Hyp4a} we have that $e^{\Lambda^{a}(x,t)}>0$ for any $x\in E$, $a\in \mathbb{U}(x)$, $0\leq t < t_{*}(x)$
($0\leq t \leq t_{*}(x)$ if $t_{*}(x)< \infty$). A consequence of item b) of Assumption \ref{Hyp4a} is
that for any
$x\in E$ with $t_{*}(x)=\infty$, $\lim_{t\rightarrow \infty} e^{-\alpha t-\int_0^{t}\underline{\lambda}(\phi(x,s))ds} =0$.
Therefore, for any $x\in E$ with $t_{*}(x)=\infty$, $A\in \mathcal{B}(E)$, $\widehat{a}=(a,a_\partial) \in \mathbb{S}(x)$,
$w\in \mathbb{M}(\partial E\times \mathbb{U})_{+}$, we have that
$G(x,\widehat{a};A) = \int_0^{t_{*}(x)}e^{-\alpha s - \Lambda^{a}(x,s)}\lambda QI_{A}(\phi(x,s),\ell(x,a,s)) ds$, and
$Hw(x,\widehat{a}) =0$.
\end{remark}

\begin{definition}\label{Vcal} For any $\mu\in \mathcal{S}_{\mathbb{S}^{r}}$ the kernel on $E$ given $E$ is defined by
\begin{eqnarray*}
G_{\mu}(x;dy) & = & \int_{\mathbb{S}(x)} G(x,\widehat{a};dy) \mu(x;d\widehat{a}).
\end{eqnarray*}
For any $U=\{\eta_j;j\in \NN\}\in \mathcal{U}$ and $k\in \NN$, let us introduce the kernel $\mathcal{G}^{k}_{U}$ on $E$ given $E$ by
\begin{eqnarray*}
\mathcal{G}^{k}_{U} (x;dy) & = & {G}_{\eta_{0}} {G}_{\eta_{1}} \ldots {G}_{\eta_{k}}(x;dy)
\end{eqnarray*}
and for notational convenience, we set $\mathcal{G}^{-1}_{U}(x;dy)=I(x;dy)$.
\end{definition}

\begin{remark}\label{rem-AUX1}
Notice that for $U^{\varphi}\in \mathcal{U}_s$, we have that for any $k\in \NN$,
\begin{align*}
G^{k}_{\varphi}(x,dy) =  \mathcal{G}^{k-1}_{U^{\varphi}}(x,dy).
\end{align*}
\end{remark}

\section{Auxiliary Results}\label{AuxRes}
The following auxiliary results will be useful in the sequel, in order to re-write our continuous-time problem in a discrete-time framework, in which the stages are defined by the jump times $T_{k}$ of the PDMP. The first result
gives an interpretation of (\ref{DefLr}), (\ref{DefHr}), in terms of the jump time $T_1$.
\begin{lemma}
\label{expressL&H}
For $x\in E$, $\widehat{a}=(a,a_{\partial})\in \mathbb{S}$, $k\in \NN$ and $h\in  \mathbb{M}(E\times \mathbb{U})_{+}$
\begin{eqnarray*}
E^{U}_{(x,\widehat{a},k)} \Big[ \int_{0}^{T_{1}} e^{-\alpha s} h(\phi(x,s),a)ds \Big] & = & Lh(x,\widehat{a}), \\
E^{U}_{(x,\widehat{a},k)} \Big[ e^{-\alpha T_{1}} I_{\{T_{1}=t_{*}(x)\}} \Big] & = & e^{-\alpha t_{*}(x)-\Lambda^{a}(x,t_{*}(x))}, \\
E^{U}_{(x,\widehat{a},k)} \Big[ e^{-\alpha T_{1}} h(Z_{1}) \Big] & = & Gh(x,\widehat{a}).
\end{eqnarray*}
\end{lemma}
\noindent {\bf{Proof}:}
It is an immediate application of (\ref{distri-T}) and (\ref{distri-Q}).
\hfill$\Box$

\bigskip

The next result re-writes the cost $\mathcal{D}_i^{\alpha}(U,\nu_0)$ in a discrete-time fashion, using the operators $L$ and $H$ defined in (\ref{DefLr}) and  (\ref{DefHr}) respectively.

\begin{proposition}\label{prop-AUX1} Consider Assumptions \ref{Hyp2a}, \ref{Hyp6a}, \ref{Hyp4a}. For $U=\{\eta_k;k\in \NN\}\in\mathcal{U}$, we have that
\begin{align}
\mathcal{D}_i^{\alpha}(U,\nu_0) &= \sum_{k=0}^\infty E^{U}_{\nu_0}\Big[e^{-\alpha T_k }\Big( L f_i(Z_k,\Theta_k)+ H r_i(Z_k,\Theta_k)\Big)\Big]
\nonumber \\& = \sum_{k=0}^\infty E^{U}_{\nu_0}\Big[e^{-\alpha T_k }\Big( L f_i(Z_k,\eta_k)+ H r_i (Z_k,\eta_k)\Big)\Big]. \label{aux1}
\end{align}
\end{proposition}
\noindent {\bf{Proof}:} From the monotone convergence theorem we have that
\begin{align}
\mathcal{D}_i^{\alpha}(U,\nu_0)  = E^{U}_{\nu_0} \Bigg[ \sum_{k=1}^\infty E_{\nu_0}^{U}\Big[ &
{\int_{T_{k-1}}^{T_k} e^{-\alpha s}f_i\bigl(X(s), \theta_{N(s)} \bigr) ds}\nonumber \\
& {+}\> {\int_{]T_{k-1},T_k]}e^{-\alpha s} r_i(X(s-),\theta_{N(s-),\partial}) dp^{*}(s)| \mathcal{F}_{T_{k-1}}}\Big]\Bigg]. \label{eq1}
\end{align}
Moreover, denoting $S_{k}=T_{k}-T_{k-1}$, we have that
\begin{align*}
E_{\nu_0}^{U} & \Big[
\int_{T_{k-1}}^{T_k} e^{-\alpha s}f_i\bigl(X(s), \theta_{N(s)} \bigr) ds
+ \int_{]T_{k-1},T_k]}e^{-\alpha s} r_i(X(s-),\theta_{N(s-),\partial}) dp^{*}(s)| \mathcal{F}_{T_{k-1}} \Big] \nonumber \\
&=e^{-\alpha T_{k-1}}E_{\nu_0}^{U} \Big[ \int_{0}^{S_{k-1}} e^{-\alpha s} f_i\bigl(\phi(Z_{k-1},s), \theta_{k-1} \bigr) ds | \mathcal{F}_{T_{k-1}} \Big] \nonumber \\
&\phantom{=} + E_{\nu_0}^{U} \Big[ e^{-\alpha T_{k}} r_i \bigl(\phi(Z_{k-1},S_{k}), \theta_{k-1,\partial} \bigr) I_{\{S_{k}=t_{*}(Z_{k-1})\}}  | \mathcal{F}_{T_{k-1}} \Big] \nonumber \\
&=e^{-\alpha T_{k-1}}E_{\nu_0}^{U} \Big[ \int_{0}^{S_{k-1}} e^{-\alpha s} f_i\bigl(\phi(Z_{k-1},s), \theta_{k-1} \bigr) ds | \mathcal{F}_{T_{k-1}} \Big] \nonumber \\
&\phantom{=} + e^{-\alpha T_{k-1}} E_{\nu_0}^{U} \Big[ e^{-\alpha S_{k}} I_{\{S_{k}=t_{*}(Z_{k-1})\}}  | \mathcal{F}_{T_{k-1}} \Big]
r_i \bigl(\phi(Z_{k-1},t_{*}(Z_{k-1})), \theta_{k-1,\partial} \bigr)
\end{align*}
Now from Lemma \ref{expressL&H} and by using the Markov property of the process $\{\widehat{X}^{U}(t)\}_{t\in \RR_{+}}$, we obtain that
\begin{align}
E_{\nu_0}^{U} & \Big[
\int_{T_{k-1}}^{T_k} e^{-\alpha s}f_i\bigl(X(s), \theta_{N(s)} \bigr) ds
+ \int_{]T_{k-1},T_k]}e^{-\alpha s} r_i(X(s-),\theta_{N(s-),\partial}) dp^{*}(s)| \mathcal{F}_{T_{k-1}} \Big] \nonumber \\
&= e^{-\alpha T_{k-1}}\Big( L f_i(Z_{k-1},\Theta_{k-1})+ H r_i (Z_{k-1},\Theta_{k-1})\Big).\label{eq2}
\end{align}
Notice also that from equation (\ref{distri-Q}) we have that
\begin{align}
E_{\nu_0}^{U}\Big[e^{-\alpha T_{k-1}}\Big( L f_i(Z_{k-1}, \Theta_{k-1})+ H r_i & (Z_{k-1},\Theta_{k-1})\Big)|T_{k-1},Z_{k-1}\Big]\nonumber \\&
= e^{-\alpha T_{k-1}}\Big( L f_i(Z_{k-1},\eta_{k-1})+ H r_i (Z_{k-1},\eta_{k-1})\Big).\label{eq3}
\end{align}
Combining (\ref{eq1}) with (\ref{eq2}) and (\ref{eq3}) we get the desired result.
\hfill$\Box$

\bigskip

The next result establishes a connection between the operator $\mathcal{G}^{k-1}_{U}$ presented in Definition \ref{Vcal} and $E_{\nu_0}^U \big[ e^{-\alpha T_k} h(Z_k) \big]$.

\begin{proposition}\label{prop-aux-1aa} Consider Assumptions \ref{Hyp2a}, \ref{Hyp6a}, \ref{Hyp4a}.
For any $U=\{\eta_j;j\in \NN\}\in \mathcal{U}$, $h\in \mathbb{M}(E)_{+}$ and $k\in \NN$
\begin{align}\label{induction-ka}
E_{\nu_0}^U \big[ e^{-\alpha T_k} h(Z_k) \big] &= \int_{E} \mathcal{G}^{k-1}_{U}h(x)\nu_0(dx).
\end{align}
\end{proposition}
\noindent {\bf{Proof}:} For $k=0$, equation (\ref{induction-ka}) follows after noticing from (\ref{expect-nu0}) that
\begin{align*}
E_{\nu_0}^U \big[ e^{-\alpha T_0}h(Z_0) \big] = E_{\nu_0}^U \big[ h(Z_0) \big]= \int_{E} h(x)\nu_0(dx) =
\int_{E} \mathcal{G}^{-1}_{U}h(x) \nu_0(dx).
\end{align*}
Suppose (\ref{induction-ka}) holds for $k$. Then,
\begin{align*}
E_{\nu_0}^U \big[ e^{-\alpha T_{k+1}}h(Z_{k+1}) \big] & = E_{\nu_0}^U \big[ e^{-\alpha T_{k}} E_{\nu_0}^U \big[ e^{-\alpha(T_{k+1} - T_{k} )}h(Z_{k+1})| \mathcal{F}_{T_{k}} \big]\big].
\end{align*}
Now from Lemma \ref{expressL&H} and by using the Markov property of the process $\{\widehat{X}^{U}(t)\}_{t\in \RR_{+}}$, we obtain that
\begin{align*}
E_{\nu_0}^U \big[ e^{-\alpha T_{k+1}}h(Z_{k+1}) \big] &
=E_{\nu_0}^U \big[ e^{-\alpha T_{k}}Gh(Z_k,\Theta_k) \big] = E_{\nu_0}^U \big[ e^{-\alpha T_{k}} G_{\eta_{k}} h(Z_k) \big] \\
&= \int_{E} \mathcal{G}^{k-1}_{U} G_{\eta_{k}} h(x)\nu_0(dx) = \int_{E} \mathcal{G}^{k}_{U} h(x)\nu_0(dx)
\end{align*}
showing (\ref{induction-ka}) for $k+1$.
\hfill$\Box$

\bigskip
The next result combines the previous results for the stationary control case.

\begin{proposition}\label{prop-AUX1a}  Consider Assumptions \ref{Hyp2a}, \ref{Hyp6a}, \ref{Hyp4a}. For $U^\varphi\in \mathcal{U}_s$ we have that
\begin{align}
\mathcal{D}_i^{\alpha}(U^\varphi,\nu_0) &= \sum_{k=0}^\infty \int_E \Big[G^{k}_{\varphi}\Big( L f_i(.,\varphi) + H r_i (.,\varphi)\Big)\Big](x)\nu_0(dx). \label{aux1-FO}
\end{align}
\end{proposition}
\noindent {\bf{Proof}:} From Proposition \ref{prop-AUX1} we have that
\begin{align*}
\mathcal{D}^{\alpha}_i(U^\varphi,\nu_0) &= \sum_{k=0}^\infty E^{U^\varphi}_{\nu_0}\Big[e^{-\alpha T_k }\Big( L f_i(Z_k,\varphi)+ H r_i (Z_k,\varphi)\Big)\Big].
\end{align*}
From Proposition \ref{prop-aux-1aa} and Remark \ref{rem-AUX1}, for any $k\in \NN$
\begin{align*}
E^{U^\varphi}_{\nu_0}\Big[e^{-\alpha T_k }\Big( L f_i(Z_k,\varphi)+ H r_i (Z_k,\varphi)\Big)\Big]&=
\int_E \Big[G^{k}_{\varphi}\Big( L f_i(.,\varphi) + H r_i (.,\varphi)\Big)\Big](x) \nu_0(dx)
\end{align*}
completing the proof.
\hfill$\Box$

\bigskip
We define next the occupation measure for our problem.
Consider $\Gamma \in \mathcal{B}(\mathcal{K})$. We define a measure  {$\mu_{\nu_0}^U \in \mathfrak{M}(\mathcal{K})_{+}$} 
as follows:
\begin{equation}\label{measure-mu}
\mu_{\nu_0}^U(\Gamma) = \sum_{k=0}^\infty E^{U}_{\nu_0}\Big[e^{-\alpha T_k } I_{\Gamma}(Z_k,\Theta_k)\Big].
\end{equation}
For any $\mu \in \mathfrak{M}(\mathcal{K})_{+}$  we denote, for notational convenience, $\widetilde{\mu}$ as the marginal of $\mu$ on $E$.

\begin{proposition}  Consider Assumptions \ref{Hyp2a}, \ref{Hyp6a}, \ref{Hyp4a}.
For any $B\in \mathcal{B}(E)$, we have that
\begin{equation}\label{aux3}
\widetilde{\mu}_{\nu_0}^U(B) = \nu_0(B) + \int_{E\times \mathbb{S}} G(z,\widehat{a};B) d\mu_{\nu_0}^U (z,\widehat{a}).
\end{equation}
\end{proposition}
\noindent {\bf{Proof}:} From (\ref{measure-mu}) we get that
\begin{align}
\widetilde{\mu}_{\nu_0}^U(B) & = \mu_{\nu_0}^U(B \times \mathbb{S})  = \nu_0(B) + \sum_{k=1}^\infty
E^{U}_{\nu_0}\Big[ e^{-\alpha T_k} I_{B\times \mathbb{S}}(Z_k,\Theta_k)\Big] \nonumber \\
& = \nu_0(B) + \sum_{k=1}^\infty  E^{U}_{\nu_0}\Big[ e^{-\alpha T_{k-1}} E_{\nu_0}^{U}\Big[e^{-\alpha(T_k-T_{k-1})}
I_{B\times \mathbb{S}}(Z_k,\Theta_k)|\mathcal{F}_{T_{k-1}}\Big]\Big].\label{aux4}
\end{align}
Notice now that
\begin{align}\label{aux5}
E_{\nu_0}^{U}\Big[e^{-\alpha(T_k-T_{k-1})}I_{B\times \mathbb{S}}(Z_k,\Theta_k)| \mathcal{F}_{T_{k-1}} \Big]  & =
E_{\nu_0}^{U}\Big[e^{-\alpha(T_k-T_{k-1})}I_{B}(Z_k)|\mathcal{F}_{T_{k-1}} \Big] \nonumber \\
&= G(Z_{k-1},\Theta_{k-1};B).
\end{align}
Combining (\ref{aux4}) and (\ref{aux5}) 
we get that
\begin{align*}
\widetilde{\mu}_{\nu_0}^U(B) & = \nu_0(B) + \sum_{k=1}^\infty  E^{U}_{\nu_0}\Big[ e^{-\alpha T_{k-1}} G(Z_{k-1},\Theta_{k-1};B)\Big] \nonumber \\
&= \nu_0(B) + \sum_{k=0}^\infty  E^{U}_{\nu_0}\Big[ e^{-\alpha T_{k}} G(Z_{k},\Theta_{k};B)\Big]\nonumber \\
&= \nu_0(B) + \int_{E\times \mathbb{S}} G(z,\widehat{a};B) d\mu_{\nu_0}^U (z,\widehat{a}),
\end{align*}
completing the proof.
\hfill$\Box$

\section{Equivalence Between the Constrained and the Linear Programming Problems}\label{sec6}

In this section we introduce the LP formulation, presented in (\ref{eq-objP})-(\ref{eq-resP1}), and the first main result of the paper, Theorem \ref{consistency},  which establishes the equivalence between the original problem and the LP problem. Define the functions $w$ on $\mathcal{K}$ and $w_0$ in $E$ as follows: for $(x,\widehat{a}) \in \mathcal{K}$, and arbitrary $c_0>0$,
\begin{align}
w(x,\widehat{a}) &\doteq c_0 + L f_0(x,\widehat{a})+ H r_0 (x,\widehat{a})>0,\label{defw}\\
w_0(x) & \doteq \inf_{\widehat{a} \in  \mathbb{S}(x)} w(x,\widehat{a}) = c_0 + \inf_{\widehat{a}  \in \mathbb{S}(x)}
\Big\{ L f_0(x,\widehat{a})+ H r_0 (x,\widehat{a})\Big\} >0.
\label{defw0}
\end{align}

We will need the following assumption, in which item b) is somehow related to the so-called expected growth condition
(see, for instance, Assumption 3.1 in \cite{guo06} for the discrete-time case, or
Assumption A in \cite{GUO-2011} for the continuous-time case).
\begin{assumption}
\label{asum-sum}
\begin{itemize}
\item[a)] The mappings $w$ and $w_{0}$ satisfy $w\in \mathbb{M}(\mathcal{K})_{+}$ and $w_{0}\in \mathbb{M}(E)_{+}$.
\item[b)]There exist $b\in \mathbb{M}(E)$,
$c> -\alpha$, and $v\in\mathbb{M}^{ac}(E)\inter \mathbb{B}_{w_0}(E)_{+}$, such that for any $(x,\widehat{a})\in \mathcal{K}$
with $\widehat{a}=(a,a_\partial)$, $Lb(x,\widehat{a})$ is well defined with values in $\mathbb{R}$, and the following inequalities are satisfied:
\begin{align}
\mathcal{X}v( & \phi(x,t)) + cv(\phi(x,t)) \nonumber \\
& -\lambda(\phi(x,t),\ell(x,t,a))\Big[v(\phi(x,t)) - Qv(\phi(x,t),\ell(x,t,a))\Big] \leq b(\phi(x,t))
\label{eqp1=},
\end{align}
\begin{align}
&\lambda(\phi(x,t),\ell(x,t,a))  + \frac{1}{c+\alpha}b(\phi(x,t)) \leq v(\phi(x,t)), \label{eqp1a=}
\end{align}
for $t\in[0,t_{*}(x))$ and
\begin{eqnarray}
v(\phi(x,t_{*}(x))) \geq Q v(\phi(x,t_{*}(x)),a_{\partial}) + c+\alpha,
\label{eqp2=}
\end{eqnarray}
for $t_{*}(x)<\infty$.
\end{itemize}
\end{assumption}

\begin{remark}\label{rem-nuw0} From Proposition \ref{prop-AUX1}, we have that for any $U\in \mathcal{U}$,
\begin{align}
\int_E w_0(x) d\nu_0(x) & \leq  c_0 + \int_E \inf_{\widehat{a}  \in \mathbb{S}(x)} \Big\{ L f_0(x,\widehat{a})+
H r_0 (x,\widehat{a})\Big\} d\nu_0(x)\nonumber \\& \leq  c_0 + \mathcal{D}^\alpha_0(U,\nu_0)
\leq c_0 + \mathcal{J}_{\mathcal{D}}^{\alpha}(\nu_0).\label{nuw}
\end{align}
From Assumption \ref{Hypdis} and (\ref{nuw}) we get $\int_E w_0(x) d\nu_0(x)<\infty$, that is, $\nu_0 \in \mathcal{P}_{w_0}(E)$.
\end{remark}

\begin{remark}\label{rem-nuw1} Notice that for $v$ as in item b) of Assumption \ref{asum-sum} and any $\nu\in \mathcal{P}_{w_0}(E)$, we have that $\nu(v) = \int_E v(x) \nu(dx)\leq \|v\|_{w_0} \int_E w_0(x) \nu(dx)<\infty$.
In particular, from Remark \ref{rem-nuw0}, $\nu_0(v)<\infty$.
\end{remark}

\begin{remark}\label{mu-mutil} If $\mu \in \mathcal{M}_w(\mathcal{K})_{+}$ then $\widetilde{\mu} \in \mathcal{M}_{w_0}(\mathcal{K})_{+}$ since
\begin{equation*}
\int_E w_0(x) d\widetilde{\mu}(x) \leq \int_\mathcal{K}  [c_0+w(x,\widehat{a})] d\mu(x,\widehat{a}) < \infty.
\end{equation*}
\end{remark}

We introduce now a linear programming formulation for the constrained problem posed in (\ref{al-disc1}).
\begin{definition}
The Problem {\bf{P}} is defined as follows:
\begin{align}
\inf_{\mu\in \mathbf{L}} \int_\mathcal{K}\Big( L f_0(x,\widehat{a}) + H r_0(x,\widehat{a}) \Big) d\mu(x,\widehat{a})\label{eq-objP}
\end{align}
where $\mathbf{L}$ is defined as the set of measure $\mu\in \mathcal{M}_w(\mathcal{K})_{+}$ satisfying for any $B\in\mathcal{B}(E)$
\begin{align}
\widetilde{\mu}(B) - \int_\mathcal{K} G(x,\widehat{a};B) d\mu(x,\widehat{a}) = \nu_0(B),\label{eq-resP}
\end{align}
and
\begin{align}
\int_\mathcal{K}\Big( L f_i(x,\widehat{a}) + H r_i(x,\widehat{a}) \Big) d\mu(x,\widehat{a})\leq d_i,\,\,\,i\in \NN_{n}^{*},
\label{eq-resP1}
\end{align}
\end{definition}

In what follows we set $\mathcal{L}(x,\widehat{a}) \doteq LI_{E\times\widetilde{\mathbb{U}}}(x,\widehat{a})$ and
$\mathcal{H}(x,\widehat{a}) \doteq HI_{E\times\widetilde{\mathbb{U}}}(x,\widehat{a})$ for all $(x,\widehat{a})\in \mathcal{K}$.
From Assumption \ref{Hyp4a} (b) it follows that $ \mathcal{L}(x,\widehat{a}) \leq K_\lambda$.
Notice that we have the following identities for any $\eta_0 \in \mathcal{S}_{\mathbb{S}^{r}}$:
\begin{align}
I_E(x) &= L(\lambda + \alpha)(x,\eta_0) + \mathcal{H}(x,\eta_0),\label{identity1}\\
I_E(x)&=G(x,\eta_0;E) + \alpha \mathcal{L}(x,\eta_0),\label{identity2}
\end{align}
for any $x\in E$.
The following result was proved in \cite{book-Springer}, using an hypothesis similar to Assumption \ref{asum-sum}:

\begin{proposition} Consider Assumptions \ref{Hyp2a}, \ref{Hyp6a}, \ref{Hyp4a}, \ref{asum-sum}. For all $(x,\widehat{a})\in \mathcal{K}$,
\begin{equation}\label{ineqvx}
v(x) \geq L((c+\alpha)v-b)(x,\widehat{a}) + (c+\alpha)\mathcal{H}(x,\widehat{a})+Gv(x,\widehat{a}).
\end{equation}
\end{proposition}
\noindent {\bf{Proof}:}
See Proposition 4.26 in \cite{book-Springer}.
\hfill$\Box$

\bigskip

We have the following proposition:

\begin{proposition}\label{prop-qU} Consider Assumptions \ref{Hyp2a}, \ref{Hyp6a}, \ref{Hyp4a}, \ref{asum-sum}.
For any $U=\{\eta_k;k\in \NN\}\in\mathcal{U}$ and $x \in E$, we have that
\begin{align}
\label{Sum-Gcal}
0 \leq \sum_{j=-1}^{\infty}\mathcal{G}^{j}_{U}(x,E) \leq \frac{1}{c+\alpha} v(x) + 1,
\end{align}
and, as a consequence, it follows that
\begin{align}
\sum_{k=0}^\infty E_{\nu_0}^U \big[e^{-\alpha T_k} \big] & \leq  \frac{1}{c+\alpha} \nu_0(v) + 1.
\label{induction-kb}
\end{align}
\end{proposition}
\noindent {\bf{Proof}:}
For $U=\{\eta_k;k\in \NN\}\in\mathcal{U}$ define the following sequence:
\begin{align}
q^U_{k}(x) &= \sum_{j=-1}^{k-1}\mathcal{G}^{j}_{U}(x,E),
\label{qkU}
\end{align}
for $k\in \NN$.
Notice that from Definition \ref{Vcal} and (\ref{qkU}) we have that for $k\in \NN$
\begin{align}
q^U_{k+1}(x) = G_{\eta_{0}} q^{U^\prime}_k(x) + I_E(x)\label{qkUprime}
\end{align}
where $U^\prime = \{\eta_k^\prime;k\in \NN\}$ with $\eta_k^\prime = \eta_{k+1}$.
Let us show by induction that
\begin{align}
\label{ineq4}
0 \leq q^{U}_{k}(x) \leq \frac{1}{c+\alpha} v(x) + 1.
\end{align}
For $k=0$ we have from (\ref{qkU}) that $q^{U}_{0}(x) = I_E(x) \leq \frac{1}{c+\alpha} v(x) + 1$ since by assumption $v$ is positive. Suppose (\ref{ineq4}) holds for $k$. From (\ref{eqp1a=}) and (\ref{ineqvx}) we have that $0\leq L((c+\alpha)v-b)(x,\widehat{a})\leq v(x)$ and thus $0\leq L((c+\alpha)v-b)(x,\eta_0)\leq v(x)$. Notice also that $0\leq \mathcal{H}(x,\eta_0)\leq 1$. From (\ref{ineqvx}) we get that
\begin{equation}\label{ineqvx1}
G_{\eta_{0}}v(x) \leq v(x)+ L(b - (c+\alpha)v)(x,\eta_0) - (c+\alpha)\mathcal{H}(x,\eta_0).
\end{equation}
From the induction hypothesis we have that $q^{U^\prime}_k(x) \leq \frac{1}{c+\alpha} v(x) + 1$.
Combining equations (\ref{identity1}), (\ref{qkUprime}) and (\ref{ineq4}), we obtain that
\begin{align*}
q^U_{k+1}(x) &\leq G_{\eta_{0}}\Big( \frac{1}{c+\alpha} v +  1\Big)(x) + L(\lambda + \alpha)(x,\eta_0) + \mathcal{H}(x,\eta_0).
\end{align*}
From equation (\ref{ineqvx1}),  we get that
\begin{align*}
q^U_{k+1}(x) & \leq \frac{1}{c+\alpha} \Big[ v(x)+ L(b - (c+\alpha)v)(x,\eta_0) - (c+\alpha)\mathcal{H}(x,\eta_0) \Big] \\
& \quad +  G_{\eta_{0}}(x,E) + L(\lambda + \alpha)(x,\eta_0) + \mathcal{H}(x,\eta_0)\\
& = \frac{1}{c+\alpha}v(x) - L\Big(v-\lambda-\frac{b}{c+\alpha}\Big)(x,\eta_0) + G_{\eta_{0}}(x,E)+\alpha \mathcal{L}(x,\eta_0).
\end{align*}
Now, observe that from (\ref{eqp1a=}) we obtain $L\Big(v-\lambda-\frac{b}{c+\alpha}\Big)(x,\eta_0)\geq 0$ and that from (\ref{identity2}),
$G_{\eta_{0}}(x,E) +  \alpha \mathcal{L}(x,\eta_0)=1$. Consequently, we have shown  equation (\ref{ineq4}).
Now, we get (\ref{Sum-Gcal}) since $\ds 0\leq \lim_{k\rightarrow \infty} q^U_{k}(x) = \sum_{j=-1}^{\infty}\mathcal{G}^{j}_{U}(x,E) \leq \frac{1}{c+\alpha} v(x) + 1$.
Moreover, by taking $h(x) = I_E(x)$ in Proposition \ref{prop-aux-1aa} we get that for every $k\in \NN$,
\begin{equation}\label{Tkfinite0}
E_{\nu_0}^U \big[ e^{-\alpha T_k} \big] = \int_{E} \mathcal{G}^{k-1}_{U}(x,E) \nu_0(dx).
\end{equation}
completing the proof.
\hfill$\Box$

\begin{remark}
\label{RemJump3}
Notice that from Assumption \ref{asum-sum} we have that there is no accumulation point of the jump times (which is considered to be an assumption  in  \cite{davis93}) since
\begin{align*}
E_{\nu_0}^U \big[ e^{-\alpha T_k} \big] & = E_{\nu_0}^U \big[ e^{-\alpha T_k}I_{\{T_{k}\leq t\}} + e^{-\alpha T_k}I_{\{T_{k}> t\}} \big]
\geq  e^{-\alpha t}E_{\nu_0}^U \big[ I_{\{T_{k}\leq t\}} \big]
\end{align*}
so that from (\ref{induction-kb}) and the monotone convergence theorem we get that
\begin{align*}
e^{-\alpha t}E^{U}_{\nu_0}\Big[ \sum_{k=1}^{\infty} I_{\{T_{k}\leq t\}} \Big] =
e^{-\alpha t} \sum_{k=1}^{\infty} E^{U}_{\nu_0} \Big[ I_{\{T_{k}\leq t\}} \Big] \leq
\sum_{k=1}^\infty E_{\nu}^U [ e^{-\alpha T_k} ] < \infty
\end{align*}
that is, $E^{U}_{\nu_0}\Big[ \sum_{k=1}^{\infty} I_{\{T_{k}\leq t\}} \Big] < \infty$. In particular we have that
$T_{k}\rightarrow \infty$ as $k \rightarrow \infty$, $P^U_{\nu_0}$-a.s. for all $U\in \mathcal{U}$.
\end{remark}

Set $\mathcal{U}_f \doteq \{U\in\mathcal{U}: \mathcal{D}_0^\alpha(U,\nu_0)<\infty,\,\,\,\mathcal{D}_i^\alpha(U,\nu_0)\leq d_i,i\in \NN_{n}^{*}\}$.
The following theorem presents the equivalence between Problem {\textbf{P}} given by the linear programming formulation posed in (\ref{eq-objP})-(\ref{eq-resP1}), and the constrained discounted piecewise deterministic Markov process problem.
\begin{theorem}\label{consistency}
Consider Assumptions \ref{Hyp2a}, \ref{Hyp6a}, \ref{Hypdis}, \ref{Hyp4a}, \ref{asum-sum}.
We have that:
\begin{enumerate}

\item [a)]For any $U \in \mathcal{U}_f $ the measure $\mu^U_{\nu_{0}}$  defined as in (\ref{measure-mu}) is in $ \mathcal{M}_w(\mathcal{K})_{+}$ and is feasible for Problem {\textbf{P}}.
Moreover,
\begin{align}\label{equality1a}
\mathcal{D}_0^\alpha(U,\nu_0) & = \int_\mathcal{K} (L f_0(x,\widehat{a}) + H r_0(x,\widehat{a})) {\mu^U_{\nu_{0}}}(d(x,\widehat{a})) \nonumber \\
& {\geq \inf_{\mu\in \mathbf{L}} \int_\mathcal{K}\Big( L f_0(x,\widehat{a}) + H r_0(x,\widehat{a}) \Big) d\mu(x,\widehat{a}).}
\end{align}

\item [b)]  For any measure $\mu \in \mathcal{M}_w(\mathcal{K})_{+}$ feasible for Problem {\textbf{P}} there exists an admissible randomized stationary control strategy $U^\varphi \in\mathcal{U}_s \cap \mathcal{U}$ for some $\varphi \in \mathcal{S}_{\mathbb{S}^r}$ 
    such that
\begin{align*}
\mathcal{J}_{\mathcal{D}}^\alpha (\nu_0)\leq \mathcal{D}_0^\alpha(U^\varphi,\nu_0) & =  \int_\mathcal{K} (L f_0(x,\widehat{a}) + H r_0(x,\widehat{a})) \mu(d(x,\widehat{a}))\\
& =  \int_E (L f_0(x,\varphi) + H r_0(x,\varphi)) \widetilde{\mu}(dx),\\
\mathcal{D}_i^\alpha(U^\varphi,\nu_0)& \leq d_i,\,\,i\in \NN_{n}^{*}.
\end{align*}
\end{enumerate}
{Moreover, the constrained discounted piecewise deterministic Markov process problem and Problem {\textbf{P}} are equivalent, that is}
\begin{equation}\label{min-PD}
\mathcal{J}_{\mathcal{D}}^\alpha (\nu_0)= \inf_{\mu\in \mathbf{L}} \int_\mathcal{K}\Big( L f_0(x,\widehat{a}) + H r_0(x,\widehat{a}) \Big) d\mu(x,\widehat{a}).
\end{equation}
\end{theorem}

\noindent {\bf{Proof}:} Regarding item a), consider the measure $\mu^U_{\nu_{0}}$ as defined in (\ref{measure-mu}).
From (\ref{aux1}) we get that $\mu^U_{\nu_{0}}$ satisfies (\ref{eq-resP1}) and that equation (\ref{equality1a}) holds.
Moreover $\mu ^{U}_{\nu_{0}} \in \mathcal{M}_w(\mathcal{K})_{+}$ is also satisfied since, recalling that $U\in \mathcal{U}_{f}$ and from Proposition \ref{prop-qU}, we have that
\begin{align*}
\int_\mathcal{K} w(x,\widehat{a}) d\mu^U_{\nu_{0}}(x,\widehat{a}) & = c_0 \mu^U(\mathcal{K}) + \int_\mathcal{K} (L f_0(x,\widehat{a}) + H r_0(x,\widehat{a})) d\mu^U_{\nu_{0}}(x,\widehat{a})  \\
&=c_0 \sum_{k=0}^\infty E_{\nu_0}^U(e^{-{\alpha T_k}}) + \mathcal{D}_0^\alpha(U,\nu_0) < \infty.
\end{align*}
For item b), consider $\mu \in \mathcal{M}_w(\mathcal{K})_{+}$ feasible for Problem {\textbf{P}}. Since $\mu$ is finite, there exists a constant $d>0$ such that $\mu(\mathcal{K}) = d$, so that $\frac{1}{d}\mu(.)$ is a probability measure concentrated on $\mathcal{K}$. From Proposition D.8(a) in \cite{hernandez96}, there exists a stochastic kernel $\varphi \in \mathcal{S}_{\mathbb{S}^{r}}$ such that
\begin{equation}\label{eqvarphi1}
\mu(B\times C) = \int_B \varphi(x;C) d\widetilde{\mu}(x),\,\, \forall B \in \mathcal{B}(E),\,\,\,\,\forall C \in \mathcal{B}(\mathbb{S}).
\end{equation}
From (\ref{eqvarphi1}) and using the notation as in (\ref{notation-hvarphi}) we have that for $i\in \NN_{n}$,
\begin{align}
\int_\mathcal{K}\Big( L f_i(x,\widehat{a}) + H r_i(x,\widehat{a}) \Big) d\mu(x,\widehat{a})& = \int_E\Big( L f_i(x,\varphi) + H r_i(x,\varphi) \Big) d\widetilde{\mu}(x), \label{funobj-1} \\
\widetilde{\mu}(B) &= \nu_0(B) + \int_E G_{\varphi}(x,E) d\widetilde{\mu}(x),\,\,\,\forall B\in\mathcal{B}(E).
\label{restrict-1}
\end{align}
Iterating (\ref{restrict-1}) we get that
\begin{align}\label{itera1}
\widetilde{\mu}(B) &= \sum_{j=0}^{m-1} \int_E G_{\varphi}^j(x,B) d\nu_0(x) + \int_E G_{\varphi}^m(x,B) d\widetilde{\mu}(x).
\end{align}
From Proposition \ref{prop-qU} we have that
\begin{equation}\label{conzero}
\sum_{j=0}^{\infty} G_{\varphi}^{j}(x,E)= \sum_{j=-1}^{\infty}\mathcal{G}^{j}_{U^{\varphi}}(I_E)(x)\leq \frac{1}{c+\alpha} v(x) + 1
\end{equation}
which implies that $G^{j}_\varphi(x,E) \rightarrow 0$ as $j\rightarrow \infty$ for each $x\in E$.
Clearly, $G^{j}_\varphi(x,E) \leq 1$ for any $j\in \NN$.
From the dominated convergence theorem we have that
\begin{align*}
0&\leq \lim_{m\rightarrow \infty}\int_E G_\varphi^m(x,B) d\widetilde{\mu}(x)\leq \lim_{m\rightarrow \infty}\int_E G_\varphi^m(x,E) d\widetilde{\mu}(x)\leq
\int_E \lim_{m\rightarrow \infty} G_\varphi^m (x,E) d\widetilde{\mu}(x) =0,
\end{align*}
and so from (\ref{itera1}) we conclude that
\begin{equation}
\widetilde{\mu}(B) = \sum_{k=0}^{\infty} \int_E G_\varphi^k(x,B) d\nu_0(x), \label{equmuB}
\end{equation}
so that from Proposition \ref{prop-AUX1a}
\begin{align*}
\int_E\Big( L f_i(x,\varphi) + H r_i(x,\varphi) \Big) d\widetilde{\mu}(x)= \sum_{k=0}^{\infty} \int_E G_\varphi^k
\Big( L f_i(.,\varphi) + H r_i(.,\varphi) \Big)(x) d\nu_0(x) = \mathcal{D}_i^{\alpha}(U^\varphi,\nu_0).
\end{align*}
Combining this with (\ref{funobj-1}) we get that
\begin{align*}
\int_\mathcal{K}\Big( L f_i(x,\widehat{a}) + H r_i(x,\widehat{a}) \Big) d\mu(x,\widehat{a})& = \mathcal{D}_i^{\alpha}(U^\varphi,\nu_0),
\end{align*}
showing item b).
\nl
From a) and b) we have (\ref{min-PD}).
\hfill$\Box$

\section{Solvability of Problem $\mathbf{P}$}\label{Solvability}

In this section we present sufficient conditions for the solvability of Problem $\mathbf{P}$ posed in (\ref{eq-objP})-(\ref{eq-resP1}).
These conditions are based on Theorem 4.1 in \cite{dufour12},
where it was considered the constrained expected total cost MDP problem, supposing that the control space didn't depend on the state variable.
Notice that in the Appendix, Theorem \ref{exist-stat} extends the results of Theorem 4.1 in \cite{dufour12} in order to consider the case in which the control space depends on the state variable, so that it can be applied to our problem. The main ideas of the proof of Theorem \ref{theoSol} below are as follows. The results in \cite{dufour12} considered the discrete-time constrained total cost MDP, therefore without any discount factor. In order to use this formulation in our problem we have to extend the state-space by setting the new state-space as $X=E\cup\{\Delta\}$, where $\Delta$ is an auxiliary state. By doing this a new Markov kernel is defined by considering $T(x,a;\{\Delta\})= 1 - G(x,a;E)$, and a new constraint is included to force that the occupation measure over $\Delta$ is zero.
With this extension we can use Theorem \ref{exist-stat} to obtain the solvability of Problem $\mathbf{P}$.
We will consider the following assumptions, similar to the ones presented in \cite{dufour12}  and \cite{Hernandez-Lerma-2000} (but not assuming the inf-compact assumption for the cost function).

\begin{assumption}\label{Hyp6bis}
The functions $Lf_i+Hr_i$ are lower semi-continuous on $\mathcal{K}$, for $i\in\NN_{n}$.
\end{assumption}

\begin{assumption}\label{Hyp7bis}
The set $\mathbb{S}(x)$ is compact for any $x\in E$ and the multifunction $\Upsilon: E\rightarrow \mathbb{S}$ defined by $\Upsilon(x)=\mathbb{S}(x)$ is upper semicontinuous.
\end{assumption}

\begin{assumption}\label{Hyp8a} $G$ is weakly continuous, that is, $Gh\in \mathbb{C}_b(\mathcal{K})$ for every $h\in \mathbb{C}_b(E)$.
\end{assumption}

We have the following theorem.

\begin{theorem}\label{theoSol}
Under the Assumptions \ref{Hyp2a}, \ref{Hyp6a}, \ref{Hypdis}, \ref{Hyp4a}, \ref{asum-sum}, \ref{Hyp6bis}, \ref{Hyp7bis}, \ref{Hyp8a}, Problem $\mathbf{P}$ posed in (\ref{eq-objP})-(\ref{eq-resP1}) is solvable.
\end{theorem}
\noindent {\bf{Proof}:} We will show that our problem can be written in the same set up as of an expected total cost Markov decision process with constraints,
so that the results of Theorem \ref{exist-stat} in the Appendix can be applied.
Let us introduce the space $X=E\union\{\Delta\}$, $A=\mathbb{S}\union\{\Delta\}$ and $A(x)=\mathbb{S}(x)$ for $x\in E$, $A(\Delta)=\{\Delta\}$
and $\mathbb{K}=\{(x,a)\in X\times A: a\in A(x)\}$.
Let us denote by $d_{1}$ (respectively $d_{2}$) the distance on the space $E$ (respectively, $\mathbb{S}$).
Without loss of generality we consider that the distance $d_{i}$ is bounded by 1 for $i=1,2$.
On the space $X$, we consider the distance $\widehat{d}_{1}$ defined by $\widehat{d}_{1}(x,y)=d_{1}(x,y)$, $\widehat{d}_{1}(x,\Delta)=2$ for any $(x,y)\in X^{2}$ and
$\widehat{d}_{1}(\Delta,\Delta)=0$. Similarly, the distance $\widehat{d}_{2}$ on $A$ is defined by $\widehat{d}_{2}(x,y)=d_{1}(x,y)$, $\widehat{d}_{2}(x,\Delta)=2$ for any $(x,y)\in \mathbb{S}^{2}$ and $\widehat{d}_{2}(\Delta,\Delta)=0$.
Clearly, it is easy to show from the assumptions that  $\mathbb{K}$ is a measurable subset of $X\times A$, $A(x)$ is compact for any $x\in X$ and
the multifunction $\Psi: X\rightarrow A$ defined by $\Psi(x)=A(x)$ is upper semicontinuous.
Let us introduce the Markov kernel $T$ on $X$ given $\mathbb{K}$ defined by
$T(x,a;\Gamma)=G(x,a;\Gamma)$ and $T(x,a;\{\Delta\})=1-G(x,a;E)$ for any $(x,a)\in \mathcal{K}$ and $\Gamma\in \mathcal{B}(E)$ and finally
$T(\Delta,\Delta;\{\Delta\})=1$.
By using the hypothesis on $G$, it follows that the kernel $T$ is weakly continuous.
Define the mapping $\mathcal{C}_{i}$ on $\mathbb{K}$ by $\mathcal{C}_{i}(x,a)=Lf_{i}(x,a)+Hr_{i}(x,a)$ for any $x\in X$ and $a\in \mathbb{S}(x)$ and
$\mathcal{C}_{i}(\Delta,\Delta)=0$ for $i\in \NN_{n}$ and $\mathcal{C}_{n+1}$ on $\mathbb{K}$ given by $\mathcal{C}_{n+1}=I_{\{(\Delta,\Delta)\}}$.
Finally, we introduce the constraints limit $R=(d_{1},\ldots,d_{n},0)$.
Clearly, the mappings $\mathcal{C}_{i}$ are lower semicontinuous on $\mathbb{K}$ for any $i\in \NN_{n+1}$.
The constrained MDP given by $\mathcal{M}=\big(X,A,(A(x))_{x\in X},T, \mathcal{C}, R\big)$ clearly satisfies the hypotheses of Theorem \ref{exist-stat}.
Therefore, there exists a Markov kernel $\varphi^{*}$ on $X$ given $\mathbb{K}$ such that
\begin{eqnarray*}
\inf_{\ds \gamma \in \mathbb{L}} \gamma(r_{0})= \mu^{\varphi^{*}}(r_{0})
\end{eqnarray*}
where
\begin{align*}
\mathbb{L} = \Big\{ \gamma \in \mathfrak{M}(\mathbb{K})_{+} : \gamma\big((\Gamma\times A) & \inter \mathbb{K}\big) =  \nu_{0}(\Gamma)
+ \gamma T(\Gamma) \text{ for any $\gamma\in \mathcal{B}(X)$ and } \nonumber \\
& \gamma(\mathcal{C}) \leq R_{j}, \text{ for } j\in\NN_{n+1}^{*} \Big\}.
\end{align*}
and
$\ds \mu^{\varphi^{*}}(\Gamma) = \sum_{t=0}^\infty \int_{\Gamma} \varphi^{*}(x;da) \nu_{0} T_{\varphi^{*}}^{t}(dx)$
for any $\Gamma\in\mathcal{B}(\mathbb{K})$
with $T_{\varphi^{*}}$ the stochastic kernel on $X$ given $X$ defined by $T_{\varphi^{*}}(x;dy)=\int_{A(x)}T(x,a;dy) \varphi^{*}(x;da)$.
However, since $\mu^{\varphi^{*}}(\mathcal{C}_{n+1})\leq 0$, we have that $\mu^{\varphi^{*}}(\{(\Delta,\Delta)\})=0$ implying that
$\ds \mu^{\varphi^{*}}(\Gamma) = \sum_{t=0}^\infty \int_{\Gamma} \varphi^{*}(x;da) \nu G_{\varphi^{*}}^{t}(dx)$
for any $\Gamma\in\mathcal{B}(\mathcal{K})$.
Now, observe that $\ds \sum_{t=0}^\infty \int_{\Gamma} \varphi^{*}(x;da) \nu G_{\varphi^{*}}^{t}(dx)=\mu^{U^{\varphi^{*}}}_{\nu_{0}}$ and
applying item $a)$ of Theorem \ref{consistency}, it follows that $\mu^{U^{\varphi^{*}}}_{\nu_{0}} \in \mathcal{M}_w(\mathcal{K})_{+}$.
Since $\mu^{\varphi^{*}}(\{(\Delta,\Delta)\})=0$, we obtain that $\mu^{U^{\varphi^{*}}}_{\nu_{0}}$ is
feasible for Problem $\mathbf{P}$. Therefore,
\begin{eqnarray*}
\inf_{\ds \gamma \in \mathbb{L}} \gamma(\mathcal{C}_{0})= \mu^{\varphi^{*}}(\mathcal{C}_{0})=\mu^{U^{\varphi^{*}}}_{\nu_{0}}(L f_0+ H r_0)
\geq \inf_{\mu\in \mathbf{L}} \int_\mathcal{K}\Big( L f_0(x,\widehat{a}) + H r_0(x,\widehat{a}) \Big) d\mu(x,\widehat{a}).
\end{eqnarray*}
Moreover, it is easy to show that
\begin{eqnarray*}
\inf_{\ds \gamma \in \mathbb{L}} \gamma(r_{0}) \leq
\inf_{\mu\in \mathbf{L}} \int_\mathcal{K}\Big( L f_0(x,\widehat{a}) + H r_0(x,\widehat{a}) \Big) d\mu(x,\widehat{a})
\end{eqnarray*}
showing the result.
\hfill$\Box$

\section{Examples}\label{Exam}

In this section we first present in subsection \ref{subsect-finite} the finite dimensional LP formulation of the problem for the case in which the control space and post-jump location set are finite. In subsection \ref{subsection-capacity} we illustrate the application of our formulation and conditions for a capacity expansion problem.

\subsection{The finite case for the control space and post-jump location}\label{subsect-finite}

Suppose that there exists a finte number of points $\{z_1\,\ldots,z_s\}\subset E$ such that after a jump the process can only move to one of these points. Moreover assume that $\mathbb{U}$ is finite, with
$\mathbb{U}=\{u_1,\ldots,u_r\}$. Denote by $\mathcal{I}_j$ and $\mathcal{I}_{j,\partial}$ the indexes such that $\mathbb{U}(z_j) = \{u_\kappa; \kappa \in \mathcal{I}_j\}$ and
$\mathbb{U}(\phi(z_j,t_*(z_j))) = \{u_\kappa; \kappa \in\mathcal{I}_{j,\partial}\}$. In this case Problem $\mathbf{P}$ can be re-written as a finite LP over $\mu_{j,\kappa,\iota}$ as follows:
\begin{align*}
(\mathbf{P})\,\,\,\min\,\,&  \sum_{j=1}^s \sum_{\kappa\in \mathcal{I}_j} \sum_{\iota \in \mathcal{I}_{j,\partial}}\Big( L f_0(z_j,u_\kappa,u_\iota)+ H r_0(z_j,u_\kappa,u_\iota) \Big) \mu_{j,\kappa,\iota}\\
\text{subject to}\,\,\,& \sum_{\kappa \in \mathcal{I}_{j}} \sum_{\iota \in \mathcal{I}_{j,\partial}}\mu_{j,\kappa,\iota} -
\sum_{p=1}^s \sum_{\kappa\in \mathcal{I}_p} \sum_{\iota \in \mathcal{I}_{p,\partial}}
G(I_{\{z_j\}})(z_p,u_\kappa,u_\iota) \mu_{p,\kappa,\iota} = \nu_{0,j},\,\,j=1,\ldots,s\\
&\sum_{j=1}^s \sum_{\kappa\in \mathcal{I}_j} \sum_{\iota \in \mathcal{I}_{j,\partial}}\Big( L f_i(z_j,u_\kappa,u_\iota)+ H r_i(z_j,u_\kappa,u_\iota) \Big) \mu_{j,\kappa,\iota}
\leq d_i,\,\,i=1,\ldots,n,\\
& \mu_{j,\kappa,\iota} \geq 0,\,\,j=1,\ldots,s,\,\,\,\kappa\in \mathcal{I}_j,\,\,\,\iota \in \mathcal{I}_{j,\partial}.
\end{align*}

\subsection{The Capacity Expansion Problem}\label{subsection-capacity}

Capacity expansion consists of general processes of adding facilities to meet, by consecutive construction of expansion projects, a rising demand.
The interested reader may consult the references \cite{davis87,luss82} for a survey on capacity expansion  including theoretical results and applications.
A point process models the arrivals of the demand with intensity $\lambda$ and at each arrival the demand increases in one unit.
The construction of a project is done at one of the possible rates $\gamma_j$ per unit of time, $j=1,\ldots,\kappa$,
and it is completed after the cumulative investment in the current project reaches a value $\tau$. Under completion the present level of demand is reduced in $\chi$ units. We will consider that $\chi=1$, $\tau$ does not depend on the present level of demand,
and that $\lambda$ is constant. We set $\gamma_0=0$ meaning that no construction is taking place.
We define the sets $\mathcal{N} \doteq \{0,\ldots,\kappa\}$ and $\mathcal{N}_j \doteq \{\iota\in \mathcal{N}; \iota \neq j\}$ for $j\in \mathcal{N}$ (that is, we exclude $j$ from $\mathcal{N}$).
The PDMP $\{X(t)\}$ takes place on $E\doteq[0,\tau)\times \NN \times \mathcal{N}$, with the vector $x=(x_{1},x_{2},x_3)\in E $ having the following interpretation: $x_1\in [0,\tau)$ denotes the amount of cumulative investment in the current project, $x_2 \in \NN$ represents the demand level and $x_3=j$ represents that construction is taking place at the rate $\gamma_j$ ($j=0$ means that no construction is taking place).

Following the notation as in sub-section \ref{prob-form-PDMP}, the parameters of the controlled PDMP are: for $x=(s,m,j)$, $m\in \NN$, $j\in \mathcal{N}_0$, $0\leq s < \frac{\tau}{\gamma_j}$, we have that
\begin{align*}
\phi(x,t)=(s+\gamma_jt,m,j),\,\,\,\,\,\,\,t_*(x)=\frac{\tau-s}{\gamma_j},
\end{align*}
and for $x=(s,m,0)$, $\phi(x,t)=(s,m,0)$, $t_*(x)=\infty$ (no construction is taking place). We set $\widetilde{\mathbb{U}}=\mathcal{N}$,
$\mathbb{U}=[0,\tau]\times \mathcal{N}$, for $x=(s,m,j)$, $j\in\mathcal{N}_0$ we set $\mathbb{U}(x)=[s,\tau]\times \mathcal{N}_j$, for $s=\tau$ (that is, $x\in \partial E$) we set $\mathbb{U}(x)=\{0\}\times\mathcal{N}$ (for simplicity in the sequel we will omit the first argument $0$ for the boundary control, since
it will play no role in the next definitions), and for $x=(s,m,0)$ we set $\mathbb{U}(x)=\{s\}\times \mathcal{N}$. The function
$\ell \in  \mathbb{M}(E\times \mathbb{U}\times \mathbb{R}_+,\widetilde{\mathbb{U}})$ is defined as follows: for $x=(s,m,j)\in E$, $a=(s_a,j_a)\in \mathbb{U}$ and $t\in \mathbb{R}_+$,
\begin{equation*}
\ell(x,a,t)\doteq
\begin{cases}
j & \text{if } s+\gamma_j t < s_a, \\
j_a & \text{if } s+\gamma_j t \geq s_a.
\end{cases}
\end{equation*}
The reasoning behind these definitions is as follows. After a jump to a point $x = (s,m,j)$ with $j\in \mathcal{N}_0$ (that is, $1\leq j\leq \kappa$, meaning that a project is undergoing at a rate $\gamma_j$) the controller will choose a time $s_a \in [s,\tau]$ and a new rate $j_a$ such that if the next
jump time due to a new arrival occurs before $s_a$ then the construction rate will be kept unchanged at $\gamma_j$, otherwise it will be changed to $\gamma_{j_a}$ (note that $\mathbb{U}(x)=[s,\tau]\times \mathcal{N}_j$, thus $j_a\neq j$). The choice $s_a=s$ means that the rate will be changed to $\gamma_{j_a}$ after the next jump time due to a new arrival, which is the case for $m=0$ (that is, $\ell(x,a,t)=j_a$ for $m=0$ since in this case $s_a=s$ is the only possible choice). The choice $s=\tau$ means that the construction rate will be kept unchanged at $\gamma_j$
after the next jump time due to a new arrival. The next construction rate after completion of a project is given by $\gamma_{a_\delta}$, with $a_\delta \in \mathcal{N}$. From these definitions we have that the control variable $\widehat{a}=(a,a_\partial)\in \mathbb{S}(x)$, $a=(s_a,j_a)$
will act on the transition measure as follows. Set for simplicity $\ell_t = \ell(x,a,t)$. For $(s,m,j) \in E$, with $j\in \mathcal{N}_0$,
\begin{align*}
&Q\bigl((s+\gamma_j t,m,j),\ell_t;A\bigr)= \delta_{\{(s+\gamma_j t,m+1,\ell_t)\}}(A) \text{ for }s\in [0,\tau),\, j\in \mathcal{N}_0,
\\
&Q\bigl((\tau,m,j),a_\partial;A\bigr)=\delta_{\{(0,m-1,a_\partial)\}}(A),\\
&Q\bigl((s,m,0),j_a;A\bigr)= \delta_{\{(s,m+1,j_a)\}}(A) \text{ for }s\in [0,\tau),\, j=0.
\end{align*}
The above equations mean that the control variable $\widehat{a}=(a,a_\partial)$ will choose the next rate as being $\gamma_{\ell_t}$ after a new demand arrives and the next rate $\gamma_{a_\partial}$ after the completion of the present project. From this we have that for $x=(s,m,j)$, $\widehat{a}=(a,a_\partial)$, $a=(s_a,j_a)$, \begin{align}
Gh(x,\widehat{a})& =\int_0^{\frac{\tau-s}{\gamma_j}} e^{-(\alpha+\lambda)}\lambda h((s+\gamma_j t,m+1,\ell_t))dt + e^{-(\alpha+\lambda)({\frac{\tau-s}{\gamma_j}})}h((0,m-1,a_\partial)), \,\,1\leq j\leq \kappa,\label{Gh1}\\
Gh(x,\widehat{a})& =\frac{\lambda}{\alpha + \lambda} h((s,m+1,j_a)),\,\,j=0.\label{Gh2}
\end{align}
The infinite-horizon expected $\alpha$-discounted costs and the constrained $\alpha$-discount value functions are as in (\ref{al-disc0}) and (\ref{al-disc1}) respectively, for positive running costs $f_i$ and boundary costs $r_i$, for $i=0,1,\ldots,n$, and a discount factor $\alpha>0$.
In this case $\alpha$ could be seen as an interest rate that brings future costs to the present value,
$\mathcal{D}_0^{\alpha}(U,\nu_0)$ could represent the total discounted cost of the project that it is desired to minimize, including penalties for
having a demand not met, while the constraint ${\mathcal{D}}_i^{\alpha}(U,\nu_0)\leq d_i$ could mean budget restrictions
that would have to be satisfied. For the initial probability measure $\nu_0$, we could consider for instance that $v_0(A)=\delta_{\{(0,0,0)\}}(A)$, meaning that initially there is no demand and no project is under construction.

It is easy to see that Assumptions \ref{Hyp2a}, \ref{Hyp4a} and \ref{Hyp7bis} are satisfied (Assumptions \ref{Hyp6a}, \ref{Hypdis}, \ref{asum-sum} a) and \ref{Hyp6bis} will depend on $f_i$ and $r_i$). Let us show next that item $b)$ of Assumption \ref{asum-sum} is satisfied. For this we consider the function
$v(s,m,j) = \lambda e^{a_1 m}$, $b(s,m,a)=0$ and $c= - \rho \alpha$ for constants $a_1>0$ and $0< \rho <1$ to be defined in the sequel.
From this definition we have that $\alpha + c= (1-\rho)\alpha > 0$ and that (\ref{eqp1a=}) is satisfied since $v(s,m,j)\geq \lambda$. It is easy
to see that (\ref{eqp1=}) and (\ref{eqp2=}) are satisfied if we have that
\begin{align}\label{eq-example}
-\rho \alpha - \lambda (1-e^{a_1}) \leq 0,\,\,\,\,
\frac{(1-\rho)\alpha}{\lambda}\leq (1-e^{-a_1}).
\end{align}
We set $a_1$ such that $e^{a_1} = 1+\frac{\alpha}{\lambda} \rho$ so that, writing $\alpha^\prime = \frac{\alpha}{\lambda}$, it follows
that (\ref{eq-example}) becomes
\begin{equation}\label{eq-example1}
1 - \rho \leq \frac{\rho}{1+\alpha^\prime \rho}, \text{ with }0<\rho<1.
\end{equation}
Defining the function $g(\rho) = \alpha^\prime \rho^2 +
(2-\alpha^\prime) \rho -1$ it follows that (\ref{eq-example1}) is satisfied if and only if $g(\rho)\geq 0$ and $0<\rho<1$. Since $g(1)=1$, $g(\frac{1}{2}) = -\frac{\alpha^\prime}{4}$ we can find $\frac{1}{2}< \bar{\rho} < 1$ such that $g(\bar{\rho})>0$, and therefore (\ref{eq-example1}) is satisfied, showing that
Assumption \ref{asum-sum} b) holds. Let us check now that Assumption \ref{Hyp8a} holds. From (\ref{Gh1}) and (\ref{Gh2}) it is easy to see that $Gh$ is bounded and
 $Gh(x_k,\widehat{a}_k) \rightarrow Gh(x,\widehat{a})$ whenever $h\in \mathbb{C}_b(E)$ and  $(x_k,\widehat{a}_k)\rightarrow (x,a)$, showing that $G$ is indeed weakly continuous.

\section*{Appendix}

In Theorem 4.1 in \cite{dufour12} it was studied the constrained total expected cost of MDPs supposing that the control space $A$ is fixed, that is,
 it doesn't depend on the state variable, while in the present paper the control space is dependent on the state variable.
The goal of this appendix is to extend in Theorem \ref{exist-stat} below the results of Theorem 4.1 in \cite{dufour12}, in order to consider the case in which the control set in the form $A(x)$, that is, depending on the state variable $x$. The basic idea will be to
 start with the control space dependent on the state-variable and, by using some extension results on the Markov kernel and cost functions, recast the case of a
 fixed control space. Moreover an extra restriction is imposed in order to obtain for the new problem occupation measures concentrated on the original control space. The final result follows after using Theorem 4.1 in \cite{dufour12}. The result in this appendix is independent of the main results in the paper and, we believe, it is interesting on its own. We will consider in this appendix a constrained MDP defined as a six-tuple $\mathcal{M}=\big(X,A,(A(x))_{x\in X},T, \mathcal{C} , R\big)$ consisting of
\begin{enumerate}
\item [(a)] a Borel space $X$ which is the state space,
\item [(b)] a Borel space $A$, representing the control or action set.
A family $\{A(x):x\in X\}$ of nonempty measurable subsets of $A$, where $A(x)$ is the set of feasible controls or actions when the system is in state $x\in X$.
We suppose that $\mathbb{K}=\{(x,a)\in X\times A: a\in A(x)\}$ is a measurable subset of $X \times A$.
\item [(c)] a stochastic kernel $T$ on $X$ given $\mathbb{K}$ which stands for the transition law of the controlled process,
\item [(d)] a measurable function $\mathcal{C}_{0}: \mathbb{K} \rightarrow \mathbb{R}$ representing the running cost,
\item [(f)] measurable functions $\mathcal{C}_{i}: \mathbb{K} \rightarrow \mathbb{R}$ for $i\in \NN_{q}^{*}$ representing the constraints,
\item[(g)] constraint limits $R=(R_{1},\ldots,R_{q})\in \RR^{q}$
\end{enumerate}
Define $H_{0}=X$ and $H_{t}=\mathbb{K}\times H_{t-1}$ for $t\geq 1$. A control policy is a sequence $\pi=(\pi_{t})_{t\in \NN}$ of stochastic kernels
$\pi_{t}$ on $A$ given $H_{t}$.
Let $\Pi$ be the class of all policies. Following standard arguments (see for example \cite[Chapter 2]{hernandez96}),
for any policy $\pi \in \Pi$ and any initial distribution $\nu$ on $X$, it can be defined a probability on the canonical space $\Omega=(X\times A)^{\infty}$, labeled $P^{\pi}_{\nu}$, and a stochastic process
$\big((x_{t},a_{t})\big)_{t\in \NN}$ where $(x_{t})_{t\in \NN}$ is the state process and $(a_{t})_{t\in \NN}$
is the control process satisfying for any $B\in \mathcal{B}(X)$, $C\in \mathcal{B}(A)$ and $t\in \NN$,
$P^{\pi}_{\nu}(x_{0}\in B)=\nu(B)$, $P^{\pi}_{\nu}(a_{t}\in C| h_{t})=\pi_{t}(h_{t};C)$, and
$P^{\pi}_{\nu}(x_{t+1}\in B|h_{t},a_{t})=T(x_{t},a_{t};B)$
where $h_{t}=(x_{0},a_{0},\ldots,x_{t-1},a_{t-1},x_{t})$.
The expectation with respect to  $P^{\pi}_{\nu}$ is denoted by $E^{\pi}_{\nu}$.
Suppose that we are given an initial distribution $\nu$ on $X$. The optimization problem we consider consists in minimizing the cost function
\begin{eqnarray}
\label{opti-1}
v(\nu,\pi)=E_{\nu}^\pi \Big[ \sum_{t=0}^\infty \mathcal{C}_{0}(x_t,a_t)\Big],
\end{eqnarray}
over the set of feasible control policies, labeled
$\Pi_{c}$, defined by the set of policies $\pi \in \Pi$ such that
\begin{eqnarray}
\label{opti-2}
v_{i}(\nu,\pi)=E_{\nu}^\pi \Big[ \sum_{t=0}^\infty \mathcal{C}_{i}(x_t,a_t)\Big] \leq R_{i},
\end{eqnarray}
for $i\in \NN_{q}^{*}$.
For a policy $\pi\in \Pi$, let us introduce the following expected state-action frequency or occupation measure induced by $\pi\in \Pi$
$\mu^\pi(\Gamma) = \sum_{t=0}^\infty P^\pi_\nu \big( (x_{t},a_{t})\in\Gamma \big)$
for any $\Gamma\in\mathcal{B}(\mathbb{K})$.

\begin{assumption}
\label{Ass1}
We assume that:
\begin{itemize}
\item[i)] The control space $A(x)$ is compact for any $x\in X$ and the multifunction
$\Psi: X \rightarrow A$ defined by $\Psi(x)=A(x)$ is upper semicontinuous.
\item[ii)] The mappings $\mathcal{C}_{i}$ for all $i\in \NN_{q}$ are non-negative and lower semi-continuous on $\mathbb{K}$.
\item[iii)] The kernel $T$ is weakly continuous, that is $Tf$ is continuous on $\mathbb{K}$ for any $f\in \mathbb{C}_{b}(X)$.
\end{itemize}
\end{assumption}

\begin{theorem}
\label{exist-stat} Suppose Assumptions \ref{Ass1} holds and that there exists a measure $\beta\in \mathbb{L}$ such that $\beta(r_{0})<\infty$ where
\begin{eqnarray*}
\mathbb{L} = \Big\{ \gamma \in \mathfrak{M}(\mathbb{K})_{+} : & \gamma\big((\Gamma\times A) \inter \mathbb{K}\big) =  \nu(\Gamma)
+ \gamma T(\Gamma) \text{ for any $\gamma\in \mathcal{B}(X)$ and } \nonumber \\
& \gamma(\mathcal{C}_{k}) \leq R_{k}, \text{ for } k\in\NN_{q}^{*} \Big\}.
\end{eqnarray*}
Then there exists a randomized stationary policy $\varphi^{*}\in \Pi_{c}$
such that
\begin{eqnarray*}
\inf_{\ds \gamma \in \mathbb{L}} \gamma(\mathcal{C}_{0})= \mu^{\varphi^{*}}(\mathcal{C}_{0}) = \inf_{\pi\in \Pi_{c}} v(\nu,\pi)  = v(\nu,\varphi^{*}).
\end{eqnarray*}
\end{theorem}
\textbf{Proof:} Without loss of generality, the set $A$ can be considered as a measurable subset of a compact Polish space $\widehat{A}$ (see for example Proposition 7 in \cite{yushkevich97}). Now by using Tietze's Theorem, the mappings $\mathcal{C}_{i}$ for all $i\in \NN_{q}$ can be extended to
non-negative, lower semi-continuous mappings defined on $X\times \widehat{A}$. The corresponding extensions will be denoted by
$\widehat{\mathcal{C}}_{i}$ for all $i\in \NN_{q}$.
The space $\mathcal{M}(X)$ is endowed with the weak topology. It is a locally convex topological vector space.
From item $iii)$ of Assumption \ref{Ass1}, the Markov kernel $T$ defines a continuous
mapping from $\mathbb{K}$ to the convex set $\mathcal{P}(X) \subset \mathcal{M}(X)$.
By using item $i)$ of Assumption \ref{Ass1}, the set $\mathbb{K}$ is closed.
Consequently, from Dugundji's Theorem (see Theorem 4.1 in \cite{dugundji51}), the mapping $T$ can be extended to a continuous mapping defined on
$X\times \widehat{A}$ and denoted by $\widehat{T}$.
Now, consider an additional constraint $\widehat{\mathcal{C}}_{q+1}$ defined on $X\times \widehat{A}$ by $\widehat{\mathcal{C}}_{q+1}=I_{\mathbb{K}^{c}}$.
It is a lower semi-continuous mapping.
Let us define the vector of constraints by $\widehat{\mathcal{C}}=(\widehat{\mathcal{C}}_{j})_{j\in\NN_{q+1}}$ and the constraint limits by $\widehat{R}=(R_{1},\ldots,R_{q},0)$.
Consider the constrained MDP defined by $\widehat{\mathcal{M}}=\big(X,\widehat{A},\widehat{T}, \widehat{\mathcal{C}},\widehat{R} \big)$.
Define $\widehat{H}_{0}=X$ and $\widehat{H}_{t}=X\times \widehat{A}\times H_{t-1}$ for $t\geq 1$. A control policy is a sequence
$\widehat{\pi}=(\widehat{\pi}_{t})_{t\in \NN}$ of stochastic kernels $\pi_{t}$ on $\widehat{A}$ given $\widehat{H}_{t}$.
Let $\widehat{\Pi}$ be the class of all policies for this model. By using the same arguments as before,
for any policy $\widehat{\pi} \in \widehat{\Pi}$, it can be defined a probability on the canonical space
$\widehat{\Omega}=(X\times \widehat{A})^{\infty}$, labeled $\widehat{P}^{\widehat{\pi}}_{\nu}$, and a stochastic process
$\big((\widehat{x}_{t},\widehat{a}_{t})\big)_{t\in \NN}$ satisfying for any $B\in \mathcal{B}(X)$, $\widehat{C}\in \mathcal{B}(\widehat{A})$ and $t\in \NN$,
$\widehat{P}^{\widehat{\pi}}_{\nu}(\widehat{x}_{0}\in B)=\nu(B)$,
$\widehat{P}^{\widehat{\pi}}_{\nu}(\widehat{a}_{t}\in \widehat{C}| \widehat{h}_{t})=\widehat{\pi}_{t}(\widehat{h}_{t};\widehat{C})$, and
$\widehat{P}^{\widehat{\pi}}_{\nu}(\widehat{x}_{t+1}\in B|\widehat{h}_{t},\widehat{a}_{t})=\widehat{T}(\widehat{x}_{t},\widehat{a}_{t};B)$
where $\widehat{h}_{t}=(\widehat{x}_{0},\widehat{a}_{0},\ldots,\widehat{x}_{t-1},\widehat{a}_{t-1},\widehat{x}_{t})$.
The expectation with respect to  $\widehat{P}^{\widehat{\pi}}_{\nu}$ is denoted by $\widehat{E}^{\widehat{\pi}}_{\nu}$.
The optimization problem we consider consists in minimizing the cost function
\begin{eqnarray}
\label{opti-hat-1}
\widehat{v}(\nu,\widehat{\pi})=\widehat{E}_{\nu}^{\widehat{\pi}} \Big[ \sum_{t=0}^\infty \widehat{\mathcal{C}}_{0}(\widehat{x}_t,\widehat{a}_t)\Big],
\end{eqnarray}
over the set of feasible control policies, labeled
$\widehat{\Pi}_{c}$, defined by the set of policies $\widehat{\pi} \in \widehat{\Pi}$ such that
\begin{eqnarray}
\label{opti-hat-2}
\widehat{v}_{i}(\nu,\widehat{\pi})=\widehat{E}_{\nu}^{\widehat{\pi}} \Big[ \sum_{t=0}^\infty \widehat{\mathcal{C}}_{i}(\widehat{x}_t,\widehat{a}_t)\Big] \leq R_{i},
\end{eqnarray}
for $i\in \NN_{q}$ and
\begin{eqnarray}
\label{opti-hat-3}
\widehat{v}_{q+1}(\nu,\widehat{\pi})=\widehat{E}_{\nu}^{\widehat{\pi}} \Big[ \sum_{t=0}^\infty I_{\mathbb{K}^{c}}(\widehat{x}_t,\widehat{a}_t)\Big] \leq 0,
\end{eqnarray}
Any measure $\gamma$ on $\mathbb{K}$ can be naturally extended to a measure on $X\times \widehat{A}$ that will be denoted by $\widehat{\gamma}$.
It can be easily shown that if $\gamma\in \mathbb{L}$ then $\widehat{\gamma}\in \widehat{\mathbb{L}}$
and $\gamma(\mathcal{C}_{0})=\widehat{\gamma}(\widehat{\mathcal{C}}_{0})$.
Conversely, for any $\widehat{\gamma}\in \widehat{\mathbb{L}}$ we have $\widehat{\gamma}(\widehat{\mathcal{C}}_{q+1})\leq 0$ and so
$\widehat{\gamma}(\mathbb{K}^{c})= 0$. Therefore, any $\widehat{\gamma}\in \widehat{\mathbb{L}}$ can be considered as a measure, labeled $\gamma$ defined on $\mathbb{K}$. It is easy to check that $\gamma\in \mathbb{L}$ and $\widehat{\gamma}(\widehat{\mathcal{C}}_{0})=\gamma(\mathcal{C}_{0})$.
Consequently,
\begin{eqnarray}
\label{inter1}
\inf_{\ds \widehat{\gamma} \in \widehat{\mathbb{L}}} \widehat{\gamma}(\widehat{\mathcal{C}}_{0})=\inf_{\ds \gamma \in \mathbb{L}} \gamma(\mathcal{C}_{0}).
\end{eqnarray}
From the previous discussion, the measure $\widehat{\beta}$ belongs to $\widehat{\mathbb{L}}$ and
$\widehat{\beta}(\widehat{\mathcal{C}}_{0})=\beta(\mathcal{C}_{0})<\infty$.
Therefore, the model $\widehat{\mathcal{M}}$ clearly satisfies the hypotheses of Theorem 4.1 in \cite{dufour12} and so,
there exists a randomized stationary policy in $\widehat{\Pi}_{c}$ defined by a Markov kernel $\widehat{\varphi}$ on $\widehat{A}$ given $X$ such that
\begin{eqnarray}
\label{inter2}
\inf_{\ds \widehat{\gamma} \in \widehat{\mathbb{L}}} \widehat{\gamma}(\widehat{\mathcal{C}}_{0})= \widehat{\mu}^{\widehat{\varphi}}(\widehat{\mathcal{C}}_{0})
\end{eqnarray}
with $\widehat{\mu}^{\widehat{\pi}}(\widehat{\Gamma}) =\sum_{t=0}^\infty \widehat{P}^{\widehat{\pi}}_\nu \big( (\widehat{x}_{t},\widehat{a}_{t})\in\widehat{\Gamma} \big),$
for any $\widehat{\Gamma}\in\mathcal{B}(X\times \widehat{A})$.
Then $\widehat{\varphi}(x;A(x))=1$ $\widehat{\mu}_{X}^{\widehat{\varphi}}-a.s.$ since $\widehat{\mu}^{\widehat{\varphi}}(\mathbb{K}^{c})=0$ and
so, there exists a Markov kernel $\varphi^{*}$ on $A$ given $X$ such that for any $x\in X$,
$\varphi^{*}(x;A(x))=1$ and $\varphi^{*}(x;\cdot)=\widehat{\varphi}(x;\cdot\inter A(x))$ $\widehat{\mu}_{X}^{\widehat{\varphi}}-a.s.$ and
$\widehat{\mu}^{\widehat{\varphi}}(\Gamma)=\mu^{\varphi^{*}}(\Gamma)$
for any $\Gamma\in \mathcal{B}(\mathbb{K})$.
Clearly, the randomized stationary policy generated by the Markov kernel $\varphi^{*}$ is in $\Pi_{c}$.
Observe, now that
\begin{eqnarray}
\label{inter3}
\inf_{\ds \gamma \in \mathbb{L}} \gamma(\mathcal{C}_{0}) \leq  \inf_{\pi\in \Pi_{c}} v(\nu,\pi)\leq v(\nu,\varphi^{*})=\mu^{\varphi^{*}}(\mathcal{C}_{0})
=\widehat{\mu}^{\widehat{\varphi}}(\widehat{\mathcal{C}}_{0})
\end{eqnarray}
Combining equations (\ref{inter1})-(\ref{inter3}), we obtain the result.

\hfill $\Box$

\bibliography{LP-PDMP}

\end{document}